%% file: KGW-ver20251031.tex
\documentclass[11pt,a4paper,reqno,dvipdfmx]{amsart}

\usepackage[truedimen,margin=25truemm]{geometry} 
\usepackage{amsmath,amssymb,amsthm,amscd}
\usepackage{extarrows,mathtools}

\title[Quantum $K$-theoretic divisor axiom for flag manifolds]
{Quantum $K$-theoretic divisor axiom for flag manifolds}
\author[C.~Lenart]{Cristian Lenart}
\address[Cristian Lenart]{Department of Mathematics and Statistics, State University of New York at Albany, 
Albany, NY 12222, U.S.A.}
\email{clenart@albany.edu}

\author[S.~Naito]{Satoshi Naito}
\address[Satoshi Naito]{Department of Mathematics, Institute of Science Tokyo, 
2-12-1 Oh-okayama, Meguro-ku, Tokyo 152-8551, Japan.}
\email{naito.s.ac@m.titech.ac.jp}

\author[D.~Sagaki]{Daisuke Sagaki}
\address[Daisuke Sagaki]{Department of Mathematics, 
Institute of Pure and Applied Sciences, University of Tsukuba, 
1-1-1 Tennodai, Tsukuba, Ibaraki 305-8571, Japan.}
\email{sagaki@math.tsukuba.ac.jp}

\author[W.~Xu]{Weihong Xu}
\address[Weihong Xu]{The Division of Physics,
Mathematics and Astronomy, California Institute of Technology, Pasadena, CA 91125, U.S.A.}
\email{weihong@caltech.edu}

\address[Leonardo C. Mihalcea]{460 McBryde Hall, Department of Mathematics, 
Virginia Tech, Blacksburg, VA 24061, USA.}
\email{lmihalce@math.vt.edu}

\makeatletter
\let\@wraptoccontribs\wraptoccontribs
\makeatother

\contrib[with an Appendix by]{Leonardo C. Mihalcea and Weihong Xu}


\allowdisplaybreaks
\numberwithin{equation}{section}

\newcommand{\Fg}{\mathfrak{g}}

\newcommand{\Ft}{\mathfrak{t}}

\newcommand{\BZ}{\mathbb{Z}}

\newcommand{\BC}{\mathbb{C}}

\newcommand{\BP}{\mathbb{P}}

\newcommand{\CO}{\mathcal{O}}
\newcommand{\CM}{\mathcal{M}}
\newcommand{\CI}{\mathcal{I}}

\newcommand{\sQ}{\mathsf{Q}}
\newcommand{\sT}{\mathsf{T}}

\newcommand{\sq}{\mathsf{q}}
\newcommand{\sB}{\mathsf{B}}

\newcommand{\vpi}{\varpi}

\newcommand{\ba}{\mathbf{a}}
\newcommand{\be}{\mathbf{e}}
\newcommand{\bp}{\mathbf{p}}
\newcommand{\bq}{\mathbf{q}}
\newcommand{\bv}{\mathbf{v}}

\newcommand{\bA}{\mathbf{A}}
\newcommand{\bB}{\mathbf{B}}

\newcommand{\bt}{\mathbf{t}}

\DeclareMathOperator{\wt}{wt}
\DeclareMathOperator{\qwt}{qwt}
\DeclareMathOperator{\ed}{end}
\DeclareMathOperator{\dec}{dec}
\DeclareMathOperator{\dist}{dist}
\DeclareMathOperator{\ev}{ev}

\newcommand{\Hom}{\mathrm{Hom}}
\newcommand{\QBG}{\mathrm{QBG}}

\newcommand{\QLS}{\mathrm{QLS}}
\newcommand{\LS}{\mathrm{LS}}
\newcommand{\Lie}{\mathrm{Lie}}

\newcommand{\edge}[1]{ \xrightarrow{\hspace{2pt}#1\hspace{2pt}} }
\newcommand{\Qe}[1]{ \xrightarrow[\mathsf{q}]{\hspace{2pt}#1\hspace{2pt}} }
\newcommand{\Be}[1]{ \xrightarrow[\mathsf{B}]{\hspace{2pt}#1\hspace{2pt}} }
\newcommand{\te}[2]{ \xrightarrow[#2]{\hspace{2pt}#1\hspace{2pt}} }

\newcommand{\lng}{w_{\circ}}
\newcommand{\mcr}[1]{\lfloor #1 \rfloor}
\newcommand{\mxcr}[1]{\lceil #1 \rceil}

\newcommand{\ha}[1]{\widehat{#1}}
\newcommand{\ol}[1]{\overline{#1}}

\newcommand{\pair}[2]{\langle #1,\,#2 \rangle}
\newcommand{\Kmet}[2]{(\!( #1,\,#2 )\!)}
\newcommand{\bra}[1]{[\![#1]\!]}

\newcommand{\dtb}[1]{\le_{#1}^{\ast}}
\newcommand{\tbmax}[3]{\max(#1W_{#2},\le_{#3}^{\ast}\nobreak)}
\newcommand{\kap}[2]{\kappa(#1,#2)}

\newcommand{\J}{J}
\newcommand{\K}{K}
\newcommand{\Lx}{L}
\newcommand{\Kc}{I \setminus \K}

\newcommand{\WJs}{W_{\J}}
\newcommand{\WLs}{W_{\Lx}}

\newcommand{\WJ}{W^{\J}}
\newcommand{\WL}{W^{\Lx}}
\newcommand{\WK}{W^{I \setminus K}}

\newcommand{\DJs}{\Delta_{\J}^{+}}
\newcommand{\DLs}{\Delta_{\Lx}^{+}}

\newcommand{\DJp}{\Delta^{+} \setminus \Delta_{\J}^{+}}
\newcommand{\DLp}{\Delta^{+} \setminus \Delta_{\Lx}^{+}}

\newcommand{\QL}{Q_{\Lx}}
\newcommand{\QJv}{Q_{\J}^{\vee}}
\newcommand{\QLv}{Q_{\Lx}^{\vee}}
\newcommand{\QKv}{Q_{\K}^{\vee}}

\newcommand{\QLvp}{Q_{\Lx}^{\vee,+}}
\newcommand{\QKvp}{Q_{\K}^{\vee,+}}

\newcommand{\bQBG}[1]{\mathbf{QBG}_{#1}^{\lhd}}
\newcommand{\bQLS}[1]{\mathbf{QLS}_{#1}^{\lhd}}
\newcommand{\pbQLS}[1]{\mathbf{QLS}_{#1}^{\K;\,\lhd}}
\newcommand{\bR}[1]{\mathbf{R}_{#1}^{\lhd}}
\newcommand{\pbR}[1]{\mathbf{R}_{#1}^{\K;\,\lhd}}
\newcommand{\bS}[1]{\mathbf{S}_{#1}^{\lhd}}
\newcommand{\pbS}[1]{\mathbf{S}_{#1}^{\K;\,\lhd}}

\newcommand{\FL}[1]{ \kappa_{\mathrm{L}}(#1) }
\newcommand{\IL}[1]{ \iota_{\mathrm{L}}(#1) }

\newcommand{\twp}[2]{ \langle \CO^{#1}, \CO_{#2} \rangle }
\newcommand{\thp}[3]{ \langle \CO^{#1}, \CO^{#2}, \CO_{#3} \rangle }

\newcommand{\dl}[1]{ \widehat{#1} }

\newcommand{\wb}{\overline}
\newcommand{\Mb}{\wb{\mathcal M}}

\newcommand{\cO}{{\mathcal O}}

\DeclareMathOperator{\pt}{pt}
\DeclareMathOperator{\Fl}{Fl}
\DeclareMathOperator{\SG}{SG}
\DeclareMathOperator{\GW}{GW}
\DeclareMathOperator{\SL}{SL}
\DeclareMathOperator{\codim}{codim}
\DeclarePairedDelimiter{\angles}{\langle}{\rangle}
\DeclarePairedDelimiter{\parens}{\lparen}{\rparen}

\DeclarePairedDelimiter\floor{\lfloor}{\rfloor}

\theoremstyle{plain}
\newtheorem{lem}{Lemma}[section]
\newtheorem{prop}[lem]{Proposition}
\newtheorem{thm}[lem]{Theorem}
\newtheorem{cor}[lem]{Corollary}
\newtheorem{conj}[lem]{Conjecture}
\newtheorem{que}[lem]{Question}

\theoremstyle{definition}
\newtheorem{dfn}[lem]{Definition}

\theoremstyle{remark}
\newtheorem{ex}[lem]{Example}
\newtheorem{rem}[lem]{Remark}

\newenvironment{enu}{%
 \begin{enumerate}%
}{\end{enumerate}}

\begin{document}


%
\begin{abstract}
We prove an identity for (torus-equivariant) 3-point, genus 0, 
$K$-theoretic Gromov-Witten invariants of flag manifolds \(G/P\), 
which can be thought of as a replacement for the ``divisor axiom'' 
in their (torus-equivariant) quantum $K$-theory. 
This identity enables us to compute these invariants 
when two insertions are Schubert classes 
and the other a Schubert divisor class.
Our type-independent proof utilizes the Chevalley formula 
for the (torus-equivariant) quantum $K$-theory ring of flag manifolds, 
which computes multiplications by Schubert divisor classes 
in terms of the quantum Bruhat graph. 
\end{abstract}

\maketitle

\section{Introduction.}

Let $G$ be a connected, simply-connected, simple (linear) algebraic group over $\BC$, 
with $T$ a maximal torus and $B$ a Borel subgroup containing $T$. 
We take an arbitrary parabolic subgroup $P$ of $G$ containing $B$, 
and let $Y \coloneqq   G/P$ be the corresponding (partial) flag manifold. 
Given an effective degree $d \in H_2(Y; \BZ)$, let $\ol{\CM}_{0,m}(Y, d)$ 
denote the Kontsevich moduli space parametrizing all $m$-point, genus $0$, 
degree $d$ stable maps to $Y$ (see \cite{FP}, \cite{Tho}). 
Given subvarieties $\Omega_1, \dots, \Omega_m \subset Y$ 
in general position, the cohomological Gromov-Witten invariant 
$\angles*{[\Omega_1], \dots, [\Omega_m]}_d^Y$ counts 
the number of parametrized curves $\BP^1 \to Y$ of degree $d$ 
with $m$ marked points in the domain (up to projective transformation) 
such that the $k$-th marked point is sent into $\Omega_k$ for $1 \leq k \leq m$,
assuming that only finitely many such ones exist. 
By transversality, one can show that if $\Omega_m = D$ is a divisor in $Y$, then 
\begin{equation}\label{eq:gwdiv}
\angles*{[\Omega_1],\dots, [\Omega_{m-1}],[D]}_d^Y = 
\left(\int_d [D]\right) \cdot \angles*{[\Omega_1],\dots, [\Omega_{m-1}]}_d^Y 
\end{equation} 
(see, for instance, \cite[Section 7]{FP}). 
This useful identity is called the divisor axiom%
\footnote{Equation \eqref{eq:gwdiv} holds for more general target spaces, 
but defining Gromov-Witten invariants in general requires the use of 
virtual fundamental classes, which we do not wish to introduce in this paper.}.
More generally, the $K$-theoretic 
Gromov-Witten (KGW) invariant 
$\angles*{[\cO_{\Omega_1}], \dots, [\cO_{\Omega_m}]}^{Y}_{d}$
is defined as the sheaf Euler characteristic 
$\chi_{\ol{\CM}_{0,m}(Y, d)}([\cO_{\GW_d}])$ of 
the Gromov-Witten subvariety $\GW_d \subset \ol{\CM}_{0,m}(Y,d)$ of stable maps that
send the $k$-th marked point into $\Omega_k$ for $1 \leq k \leq m$.
 However, no analog of the divisor axiom \eqref{eq:gwdiv} 
 is available in the general setting of $K$-theory. 

Most studies in this area have focused on $3$-point KGW invariants. 
These invariants govern the small quantum $K$-theory ring $QK(Y)$ of $Y$, 
which, introduced by \cite{Giv} and \cite{Lee}, is a deformation of 
the $K$-theory ring $K(Y)$. In particular, we have the following conjecture of 
Buch and Mihalcea (see \cite[Section~5.1.2, Conjecture~1]{Mi}), 
where $\cO^{s_i}=[\cO_{Y^{s_i}}]$ for a Schubert divisor $Y^{s_i}$ in $Y$ 
and $d_i\coloneqq  \int_d [Y^{s_i}]$.
\begin{conj}[Buch-Mihalcea] \label{conj:divisor}
When $G$ is of Lie type $A$, the following $K$-theoretic analog of the divisor axiom holds for $m=3${\rm :}
\begin{equation}\label{eq:conj}
\angles*{[\cO_{\Omega_1}], [\cO_{\Omega_2}], \cO^{s_i}}^{Y}_{d} = 
\begin{cases}
\angles*{[\cO_{\Omega_1}], [\cO_{\Omega_2}]}^{Y}_{d} & \text{\rm if } d_i > 0, \\
\angles*{\cO^{s_i}\cdot[\cO_{\Omega_1}], [\cO_{\Omega_2}]}_d^Y & \text{\rm if } d_i= 0,
\end{cases}
\end{equation} 
where $\cO^{s_i}\cdot[\cO_{\Omega_1}]$ denotes the product in the ordinary $K$-theory ring $K(Y)$.
\end{conj}

The conjecture was initially made for type $A$ flag manifolds; 
while the statement is expected to be true more generally, 
\cite[Section 4]{LM} gave a counterexample in type $G_2$ to 
the first case of \eqref{eq:conj}. 
We give details of this counterexample in Example~\ref{ex:G2-2}. 

In this paper, we prove that the second case of \eqref{eq:conj} holds true 
for all $G/P$, and the first case holds true if 
\begin{align} \label{eq:cond}
\pair{\varpi_i}{\theta^\vee} = 1,
\end{align}
where $\varpi_i$ is the fundamental weight of $G$ corresponding 
to the Schubert divisor $Y^{s_i}$, and $\theta$ is 
the highest root in the root system of $G$ (see Section~\ref{subsec:alggrp} for more details). 
This condition is satisfied for all minuscule or cominuscule fundamental weights $\vpi_i$ 
(which include all fundamental weights in type $A$), 
all fundamental weights in type $C$, and some 
but not all fundamental weights in exceptional types. 
See Section~\ref{subsec:LS} for a complete classification. 
When \eqref{eq:cond} is not satisfied, we provide an explicit formula 
for the difference between the two sides in the first case of \eqref{eq:conj}. 
The necessary corrections are given in terms of 
the quantum Lakshmibai-Seshadri (QLS) path model; 
see Section~\ref{subsec:main} for the precise statement. 
Moreover, our results hold in the more refined $T$-equivariant setting.

For $1\leq k \leq m$, let $\ev_k: \ol{\CM}_{0,m}(Y, d)\to Y$ be 
the evaluation map that sends a stable map to the image of the $k$-th marked point 
in its domain. Let $Y_u \subset Y$ be a Schubert variety 
in general position with $Y^{s_i}$, set $\cO_u \coloneqq [\cO_{Y_u}]$, and define
\begin{equation*}
\CM_d(Y_u, Y^{s_i}) \coloneqq \ev_1^{-1}(Y_u)\cap \ev_2^{-1}(Y^{s_i}) \subset \ol{\CM}_{0,3}(Y,d)
\end{equation*}
as well as
\begin{equation}\label{eq:Gamma}
\Gamma_d(Y_u, Y^{s_i}) \coloneqq   \ev_3(\CM_d(Y_u, Y^{s_i})) \subset Y.
\end{equation}
(The equivariant version of) identity~\eqref{eq:conj} is equivalent to
\begin{align}\label{eq:qclassical}
 \angles*{\cO^{s_i},\cO_u,\gamma}_d^Y=
 \chi_Y^T\parens*{\cO_{\Gamma_d(Y_u, Y^{s_i})}\cdot\gamma} \quad 
 \text{\rm for all } \gamma\in K_T(Y),
\end{align}
where $\chi_Y^T: K_T(Y)\to K_T(\pt)$ is the sheaf Euler characteristic map.
Combined with a theorem of Brion \cite{Bri}, this implies that 
when $\gamma$ comes from a suitable basis of $K_T(Y)$, 
the KGW invariant $\angles*{\cO^{s_i},\cO_u,\gamma}_d^Y$ has 
a ``positivity property''; see Corollary~\ref{cor:alt-sign} and 
Remark~\ref{rem:positivity} for more details.

One of the potential applications of 
our $K$-theoretic divisor axiom lies in proving relations 
in the ($T$-equivariant) quantum $K$-theory ring $QK_T(Y)$. 
For instance, utilizing the connection between the product $\star$ in $QK_{T}(Y)$ 
and 3-point KGW invariants, one may prove relations 
involving a product with $\cO^{s_i}$  in $QK_T(Y)$ 
from the corresponding well-understood relations in $K_T(Y)$. 
Indeed, relations in $QK_T(\SL_n(\BC)/P)$ are deduced 
in this manner in \cite{GMSXZZ}. 
Moreover, using the methods in \cite{GMSXZZ} and \cite{AHKMOX}, 
from these relations one easily obtains presentations of 
all $QK_T(\SL_n(\BC)/P)$ without relying on a known presentation of 
$QK_T(\SL_n(\BC)/B)$. Such applications may be generalizable to 
other Lie types, where presentations of the quantum $K$-theory rings are mostly unproved. 

Another potential application of the $K$-theoretic divisor axiom lies 
in deriving cancellation-free Chevalley formulas for $QK_T(G/P)$, 
since the Chevalley structure constants can be computed recursively 
from KGW invariants of the form 
$\angles*{[\cO_{\Omega_1}], [\cO_{\Omega_2}], \cO^{s_i}}^{Y}_{d}$ 
(see equation~\eqref{eq:qmulti}). This would be an alternative 
(and type-independent) approach to the one taken in~\cite{KLNS} 
for certain Grassmannians and two-step flag manifolds of type $A$.

In prior works, identity \eqref{eq:qclassical} was proved 
for cominuscule flag manifolds $G/P$ in \cite{BCMP2}, 
for the incidence variety $\Fl(1,n-1;n)$ in \cite{Xu}, 
and for $\SG(2,2n)$, the symplectic Grassmannian of lines, 
in \cite{BPX}. In each case, it was implied by the geometric statement ($\dagger$) below: 
\begin{equation} \label{eq:dagger}
\text{The general fibers of the map 
$\ev_3: \CM_d(Y_u, Y^{s_i}) \rightarrow \Gamma_d(Y_u, Y^{s_i})$ 
are rationally connected.} \tag{$\dagger$}
\end{equation}
In the cases where the correction term is nonzero 
(such as the case in Example~\ref{ex:G2-2}), 
our results imply that \eqref{eq:dagger} does not hold. 
Motivated by our findings, we would like to pose the following geometric question.
\begin{que}
  Does \eqref{eq:dagger} hold whenever $\pair{\vpi_i}{\theta^{\vee}} = 1$ or $d_i=0$?
\end{que}

In Appendix~\ref{sec:MX}, we prove the following formula, 
which may be of independent interest, 
that relates KGW invariants of $G/P$ to those of $G/B$ 
via the natural projection 
\begin{equation*}
\pi: G/B \rightarrow G/P,
\end{equation*}
generalizing a similar formula in cohomology, 
conjectured by Peterson \cite{Pet97} and proved by Woodward \cite{Woo}. 
Here, $\ha{d} \in H_2(G/B;\BZ)$ is the Peterson lift of 
$d \in H_2(G/P;\BZ)$ (see \cite[Lemma 1]{Woo} or Appendix~\ref{sec:MX} for more details).
%
%
\begin{prop} \label{prop:comparison}
For classes $\gamma_1,\dots,\gamma_m \in K_T(Y)$, we have 
\begin{equation}\label{eq:comparison}
\langle \gamma_1, \dots, \gamma_m \rangle_{d}^{G/P} = 
\langle \pi^*\gamma_1, \dots, \pi^*\gamma_m \rangle_{\dl{d}}^{G/B}.
\end{equation}
\end{prop}

Proposition~\ref{prop:comparison} reduces the problem from $G/P$ to $G/B$, 
at least in cases where the correction term vanishes 
(see Appendix~\ref{sec:MX} for more details). 
In the $G/B$ case, our main tool is the quantum $K$-Chevalley formula 
proven in \cite{NOS} and \cite{LNS} (see Theorem~\ref{thm:NOS}), 
which computes, in a cancellation-free manner, the operator 
$\cO^{s_i} \star \bullet : QK_T(G/B) \rightarrow QK_T(G/B)$ 
in the basis of Schubert classes, and the combinatorics 
is encoded in QLS paths. This formula allows us to approach the problem combinatorially, 
with the key observation that $\pair{\varpi_i}{\theta^\vee} = 1$ if and only if 
the set of QLS paths of shape $\varpi_i$ consists only of 
ordinary Lakshmibai-Seshadri (LS) paths 
(see \cite[Lemma~2.14]{MNS} and Section~\ref{subsec:LS}). 
Another tool useful for going from $G/B$ to $G/P$ is a result of Kato \cite{Kat2}, 
which states that the pushforward along $\pi$ induces a surjective ring homomorphism 
in quantum $K$-theory (see Theorem~\ref{thm:kato-morphism} for more details). 
This allows us to deduce a Chevalley formula for $QK_T(G/P)$ 
from that of $QK_T(G/B)$, although 
in the parabolic case the formula will no longer be cancellation-free.

This paper is organized as follows. 
In Section~\ref{sec:prelim}, we fix the notation for root systems, 
recall the definitions of the quantum Bruhat graph and QLS paths, 
and also give preliminaries on the quantum $K$-theory of flag manifolds. 
We state our main results in Section~\ref{sec:main}, and prove them in Section~\ref{sec:prf-main}.
Proposition~\ref{prop:comparison} and related discussions are given in Appendix~\ref{sec:MX}.

\medskip
\paragraph{\bf Acknowledgments.}
C.L. was partly supported by the NSF grants DMS-1855592 and DMS-2401755. 
S.N. was partly supported by JSPS Grant-in-Aid for Scientific Research (C) 21K03198. 
D.S. was partly supported by JSPS Grant-in-Aid for Scientific Research (C) 23K03045.
L.M. was partly supported by NSF grant DMS-2152294, and gratefully acknowledges 
the support of Charles Simonyi Endowment, which provided funding for 
the membership at the Institute of Advanced Study during 
the 2024-25 Special Year in ``Algebraic and Geometric Combinatorics''.

%
\section{Preliminaries.}
\label{sec:prelim}
%
%
\subsection{Notation for root systems.} \label{subsec:alggrp}
Let $G$ be a connected, simply-connected, simple (linear) algebraic group over $\BC$, 
$T$ a maximal torus of $G$. 
We set $\Fg \coloneqq   \Lie(G)$ and $\Ft \coloneqq   \Lie(T)$; 
$\Fg$ is a finite-dimensional simple Lie algebra over $\BC$, and 
$\Ft$ is a Cartan subalgebra of $\Fg$. 
We denote by $\pair{\cdot\,}{\cdot} : \Ft^{\ast} \times \Ft \rightarrow \BC$ 
the canonical pairing, where $\Ft^{\ast} \coloneqq   \Hom_{\BC}(\Ft, \BC)$. 
Let $\Delta \subset \Ft^{\ast}$ be the root system of $\Fg$, 
$\Delta^{+} \subset \Delta$ the set of positive roots, 
and $\{ \alpha_{j} \}_{j \in I} \subset \Delta^{+}$ the set of simple roots. 
We denote by $\alpha^{\vee} \in \Ft$ the coroot of $\alpha \in \Delta$. 
Also, we denote by $\theta \in \Delta^+$ the highest root of $\Delta$, 
and set $\rho \coloneqq   (1/2) \sum_{\alpha \in \Delta^{+}} \alpha$. 
The root lattice $Q$ and the coroot lattice $Q^{\vee}$ of $\Fg$ are defined by 
$Q \coloneqq   \sum_{j \in I} \BZ \alpha_{j}$ and $Q^{\vee} \coloneqq   \sum_{j \in I} \BZ \alpha_{j}^{\vee}$, respectively. 
We set $Q^{\vee,+}\coloneqq  \sum_{j \in I} \BZ_{\ge 0} \alpha_{j}^{\vee}$; for 
$\xi,\zeta \in Q^{\vee}$, we write $\xi \ge \zeta$ if $\xi-\zeta \in Q^{\vee,+}$. 
For $i \in I$, the weight $\vpi_{i} \in \Ft^{\ast}$ 
given by $\pair{\vpi_{i}}{\alpha_{j}^{\vee}} = \delta_{i,j}$ for all $j \in I$, 
with $\delta_{i,j}$ the Kronecker delta, 
is called the $i$-th fundamental weight. 
The (integral) weight lattice $\Lambda$ of $\Fg$ is defined by 
$\Lambda \coloneqq   \sum_{j \in I} \BZ \vpi_{j}$. 
We denote by $\BZ[\Lambda]$ the group algebra of $\Lambda$, that is, 
the associative algebra generated by formal elements $\be^{\nu}$, $\nu \in \Lambda$, 
where the product is defined by $\be^{\mu} \be^{\nu} \coloneqq   \be^{\mu + \nu}$ for $\mu, \nu \in \Lambda$. 

A reflection $s_{\alpha} \in GL(\Ft^{\ast})$, $\alpha \in \Delta$, 
is defined by $s_{\alpha} \mu \coloneqq   \mu - \pair{\mu}{\alpha^{\vee}} \alpha$ 
for $\mu \in \Ft^{\ast}$. We write $s_{j} \coloneqq   s_{\alpha_{j}}$ for $j \in I$. 
Then the (finite) Weyl group $W$ of $\Fg$ is defined to be the subgroup of $GL(\Ft^{\ast})$ 
generated by $\{ s_{j} \}_{j \in I}$, that is, $W \coloneqq   \langle s_{j} \mid j \in I \rangle$. 
For $w \in W$, there exist $j_{1}, \ldots, j_{r} \in I$ such that $w = s_{j_{1}} \cdots s_{j_{r}}$. 
If $r$ is minimal, then the product $s_{j_{1}} \cdots s_{j_{r}}$ is called a reduced expression for $w$, 
and $r$ is called the length of $w$; we denote by $\ell(w)$ the length of $w$. 

Let $\Lx$ be a subset of $I$. We set
\begin{align*}
& \QL \coloneqq   \sum_{j \in \Lx} \BZ \alpha_j, \qquad 
  \DLs \coloneqq   \Delta^{+} \cap \QL, \qquad
  \rho_{\Lx}\coloneqq  (1/2) \sum_{\alpha \in \DLs} \alpha, \\
& \QLv \coloneqq   \sum_{j \in \Lx} \BZ \alpha_j^{\vee}, \qquad
  \QLvp \coloneqq   \sum_{j \in \Lx} \BZ_{\ge 0} \alpha_j^{\vee}, \qquad 
  \WLs \coloneqq   \langle s_{j} \mid j \in \Lx \rangle. 
\end{align*}
For $w \in W$, let $\mcr{w}=\mcr{w}^{\Lx}$ and 
$\mxcr{w}=\mxcr{w}^{\Lx}$ denote the minimal(-length) coset representative
and maximal(-length) coset representative for the coset $w\WLs$, respectively. We set 
\begin{align*}
\WL \coloneqq  \bigl\{ \mcr{w}^{\Lx} \mid w \in W \bigr\} \subset W, \qquad 
\WL_{\max} \coloneqq  \bigl\{ \mxcr{w}^{\Lx} \mid w \in W \bigr\} \subset W.
\end{align*}
Also, we denote by $[\,\cdot\,]=[\,\cdot\,]_{\Lx}:
Q^{\vee} \twoheadrightarrow \QJv$ the projection 
that maps $\sum_{j \in I} c_{j}\alpha_{j}^{\vee} \in Q^{\vee}$ to 
$\sum_{j \in \Lx} c_{j}\alpha_{j}^{\vee} \in \QLv$. 

%
\subsection{The quantum Bruhat graph.}
\label{subsec:QBG}

\begin{dfn} \label{dfn:QBG} 
Let $\Lx$ be a subset of $I$. 
The (parabolic) \emph{quantum Bruhat graph} on $\WL$, 
denoted by $\QBG(\WL)$, is the ($\DLp$)-labeled
directed graph whose vertices are the elements of $\WL$ and 
whose edges are of the following form: 
$x \edge{\alpha} y$, with $x, y \in \WL$ and $\alpha \in \DLp$, 
such that $y = \mcr{x s_{\alpha}}^{\Lx}$ and either of the following holds: 
(B) $\ell(y) = \ell (x) + 1$; 
(Q) $\ell(y) = \ell (x) + 1 - 2 \pair{\rho-\rho_{\Lx}}{\alpha^{\vee}}$.
An edge satisfying (B) (resp., (Q)) is called a \emph{Bruhat edge} (resp., \emph{quantum edge}). 
When $\Lx=\emptyset$ (note that $W^{\emptyset}=W$, $\rho_{\emptyset}=0$, 
and $\mcr{x}^{\emptyset}=x$ for all $x \in W$), 
we write $\QBG(W)$ for $\QBG(W^{\emptyset})$.
\end{dfn}

For an edge $x \edge{\alpha} y$ in $\QBG(W)$, 
we sometimes write $x \Be{\alpha} y$ (resp., $x \Qe{\alpha} y$) 
to indicate that the edge is a Bruhat (resp., quantum) edge. 

Let 
\begin{equation} \label{eq:dp}
\bp: y_{0} \edge{\beta_{1}} y_{1} \edge{\beta_{2}} \cdots \edge{\beta_{r}} y_{r}
\end{equation}
be a directed path in the quantum Bruhat graph $\QBG(W)$. 
We set $\ed(\bp)\coloneqq  y_{r}$ and $\ell(\bp)=r$. 
A directed path $\bp$ is called the trivial (resp., non-trivial) one 
if $\ell(\bp)=0$ (resp., $\ell(\bp) > 0$). 
When $\bp$ is non-trivial, 
we call $\IL{\bp}\coloneqq  \beta_{1}$ and $\FL{\bp}=\beta_{r}$ 
the initial label and the final label of $\bp$, respectively. 
For $\bp$ of the form \eqref{eq:dp}, we set
\begin{equation*}
\qwt(\bp)\coloneqq  \sum_{ \begin{subarray}{c} 1 \le k \le r \\ 
\text{$y_{k-1} \edge{\beta_k} y_{k}$ is} \\
\text{a quantum edge} \end{subarray} } \beta_{k}^{\vee} \in Q^{\vee,+}. 
\end{equation*}
Let $v,w \in W$, and let $\bp$ be a shortest directed path from $v$ to $w$ in $\QBG(W)$. 
We set $\qwt(v \Rightarrow w)\coloneqq  \qwt(\bp)$; we know from \cite[Proposition~8.1]{LNSSS1} that 
$\qwt(v \Rightarrow w)$ does not depend on the choice of a shortest directed path $\bp$. 
%
%
\begin{prop}[{\cite[Proposition~8.1]{LNSSS1}}] \label{prop:qwt}
Let $v,w \in W$, and let $\bq$ be a directed path (not necessarily, shortest) 
from $v$ to $w$ in $\QBG(W)$. Then, we have $\qwt(v \Rightarrow w) \le \qwt(\bq)$.
\end{prop}

Let $\lhd$ be a reflection (convex) order on $\Delta^{+}$; 
see, e.g., \cite[Section 2.2]{KNS}.
A directed path $\bp$ of the form \eqref{eq:dp} is said to be label-increasing 
(resp., label-decreasing) with respect to $\lhd$ 
if $\beta_{1} \lhd \cdots \lhd \beta_{r}$ (resp., $\beta_{1} \rhd \cdots \rhd \beta_{r}$). 

\begin{thm}[{see, for example, \cite[Theorem~7.4]{LNSSS1}}]
Let $v,w \in W$.
\begin{enu}
\item There exists a unique label-increasing directed path from $v$ to $w$ 
in the quantum Bruhat graph $\QBG(W)$. Moreover, it is a shortest directed path 
from $v$ to $w$. 

\item There exists a unique label-decreasing directed path from $v$ to $w$ 
in the quantum Bruhat graph $\QBG(W)$. Moreover, it is a shortest directed path 
from $v$ to $w$. 
\end{enu}
\end{thm}

For $v,w \in W$, we denote by $\dec(v \Rightarrow w)$ the (unique) label-decreasing 
directed path from $v$ to $w$; note that $\qwt(v \Rightarrow w) = \qwt(\dec(v \Rightarrow w))$. 
%
%
\begin{dfn}[dual tilted Bruhat order] \label{dfn:dtilted}
For each $v \in W$, we define the dual $v$-tilted Bruhat order 
$\dtb{v}$ on $W$ as follows:
for $w_{1},w_{2} \in W$, 
%
%
\begin{equation} \label{eq:dtilted}
w_{1} \dtb{v} w_{2} \iff \ell(w_{1} \Rightarrow v) = 
 \ell(w_{1} \Rightarrow w_{2}) + \ell(w_{2} \Rightarrow v).
\end{equation}
Namely, $w_{1} \dtb{v} w_{2}$ if and only if 
there exists a shortest directed path in $\QBG(W)$ 
from $w_{1}$ to $v$ passing through $w_{2}$; 
or equivalently, if and only if the concatenation of a shortest directed path 
from $w_{1}$ to $w_{2}$ and one from $w_{2}$ to $v$ 
is one from $w_{1}$ to $v$. 
\end{dfn}
%
%
\begin{prop}[{\cite[Proposition~2.25]{NOS}}] \label{prop:tbmax}
Let $v \in W$, and let $\Lx$ be a subset of $I$. 
Then each coset $u\WL$ for $u \in W$ has a unique maximal element
with respect to $\dtb{v}$\,{\rm;} 
we denote it by $\tbmax{u}{\Lx}{v}$.
\end{prop}

The following lemma will be used in Section~\ref{subsubsec:vanish}.
%
%
\begin{lem} \label{lem:tbmax-B}
Let $\Lx$ be a subset of $I$.
Let $v,w \in W$, and assume that $\mcr{w}^{\Lx} \le \mcr{v}^{\Lx}$ in the Bruhat order. 
If $w=\tbmax{w}{\Lx}{v}$, then $w \le v$ in the Bruhat order. 
\end{lem}

\begin{proof}
%
We write $v$ as $v=\mcr{v}^{\Lx}z$ with $z \in \WLs$, and set $w'\coloneqq  \mcr{w}^{\Lx}z$. 
Since $w' \in w\WLs$ and $w=\tbmax{w}{\Lx}{v}$, 
we deduce by the definition of $\dtb{v}$ that 
$\ell(w' \Rightarrow v) = 
 \ell(w' \Rightarrow w) + \ell(w \Rightarrow v)$. 
Hence, it follows that 
\begin{equation*}
\qwt(w' \Rightarrow v) = 
\qwt(w' \Rightarrow w) + \qwt(w \Rightarrow v). 
\end{equation*}
Also, since $\mcr{w}^{\Lx} \le \mcr{v}^{\Lx}$, 
it follows that $w'=\mcr{w}^{\Lx}z \le \mcr{v}^{\Lx}z=v$, 
which implies that $\qwt(w' \Rightarrow v) = 0$. 
Since $\qwt(w' \Rightarrow w) \ge 0$ and 
$\qwt(w \Rightarrow v) \ge 0$, we deduce that 
$\qwt(w \Rightarrow v) = 0$. 
Therefore, we conclude that $w \le v$, as desired. 
This proves the lemma. 
\end{proof}

Let $\Lx$ be a subset of $I$. 
As \cite[(2.4)]{KNS}, let $\lhd$ be an arbitrary reflection 
order on $\Delta^{+}$ satisfying the condition that 
%
%
\begin{equation} \label{eq:ro}
\beta \lhd \gamma \quad 
\text{for all $\beta \in \DLs$ and $\gamma \in \DLp$}. 
\end{equation}
%
%
%
\begin{lem}[{\cite[Lemma~2.15]{KNS}}] \label{lem:tbmax}
Let $\Lx$ be a subset of $I$, and let $\lhd$ be a reflection order on $\Delta^{+}$ 
satisfying condition \eqref{eq:ro}. 
Let $v,w \in W$, and $w' \in w\WLs$. Then, 
$w' = \tbmax{w}{\Lx}{v}$ if and only if 
all the labels in the label-increasing (shortest) directed path 
from $w'$ to $v$ in $\QBG(W)$ are contained in $\DLp$. 
\end{lem}

For $w \in W$, let $\bQBG{w}$ denote 
the set of all label-increasing directed paths $\bp$ 
in $\QBG(W)$ starting at $w$, and 
satisfying the condition that 
all the labels of the edges in $\bp$ are contained in $\DLp$: 
%
%
\begin{equation} \label{eq:bQBG}
\bp: 
\underbrace{w = z_{0} \edge{\beta_{1}} z_{1} \edge{\beta_{2}} \cdots \edge{\beta_{s}} z_{s},}_{
\text{directed path in $\QBG(W)$}} \quad \text{where} \quad 
\begin{cases}
s \ge 0, \\[1mm]
\text{$\beta_{k} \in \DLp$, $1 \le k \le s$}, \\[1mm]
\beta_{1} \lhd \beta_{2} \lhd \cdots \lhd \beta_{s}. 
\end{cases}
\end{equation}
Let $\bt_{w}$ denote the trivial directed path of length $0$ 
starting at $w$ and ending at $w$; 
note that $\bt_{w} \in \bQBG{w}$. 

%
\subsection{Quantum Lakshmibai-Seshadri paths.}
\label{subsec:QLS}

Let $\lambda \in \Lambda^{+}\coloneqq  \sum_{j \in I} \BZ_{\ge 0}\vpi_{j}$ 
be a dominant (integral) weight, and take 
%
%
\begin{equation} \label{eq:J}
\J=\J_{\lambda}\coloneqq   
\bigl\{ j \in I \mid \pair{\lambda}{\alpha_{j}^{\vee}}=0 \bigr\}.
\end{equation}
%
%
\begin{dfn} \label{dfn:QBa}
For a rational number $0 \le a < 1$, 
we define $\QBG_{a\lambda}(\WJ)$ to be the subgraph of $\QBG(\WJ)$ 
with the same vertex set but having only those directed edges of 
the form $x \edge{\alpha} y$ for which 
$a\pair{\lambda}{\alpha^{\vee}} \in \BZ$ holds. 
Note that if $a = 0$, then $\QBG_{a\vpi_{i}}(\WJ) = \QBG(\WJ)$. 
\end{dfn}
%
%
\begin{dfn}[{\cite[Section~3.2]{LNSSS2}}] \label{dfn:QLS}
A \emph{quantum  Lakshmibai-Seshadri path} (QLS path for short) 
of shape $\lambda$ is a pair 
%
%
\begin{equation} \label{eq:QLS}
\eta = (\bv \,;\, \ba) = 
(v_{1},\,\dots,\,v_{s} \,;\, a_{0},\,a_{1},\,\dots,\,a_{s}), \quad s \ge 1, 
\end{equation}
of a sequence $v_{1},\,\dots,\,v_{s}$ 
of elements in $\WJ$, with $v_{k} \ne v_{k+1}$ 
for any $1 \le k \le s-1$, and an increasing sequence 
$0 = a_0 < a_1 < \cdots  < a_s =1$ of rational numbers 
satisfying the condition that there exists a directed path 
in $\QBG_{a_{k}\lambda}(\WJ)$ from $v_{k+1}$ to  $v_{k}$ 
for each $k = 1,\,2,\,\dots,\,s-1$. 
\end{dfn}

Let $\QLS(\lambda)$ denote 
the set of all QLS paths of shape $\lambda$. 
For $\eta \in \QLS(\lambda)$ of the form \eqref{eq:QLS}, 
we set $\iota(\eta)\coloneqq  v_{1} \in \WJ$, $\kappa(\eta)\coloneqq  v_{s} \in \WJ$, and 
\begin{equation} \label{eq:wt}
\wt (\eta) \coloneqq   \sum_{k=1}^{s} (a_{k}-a_{k-1}) v_{k}\lambda \in \Lambda. 
\end{equation}
Also, following \cite[(3.26) and (3.27)]{NOS}, 
for $\eta \in \QLS(\lambda)$ of the form \eqref{eq:QLS} and $v \in W$, 
we define $\kap{\eta}{v} \in W$ by the following recursive formula: 
%
%
\begin{equation} \label{eq:haw}
\begin{cases}
\ha{v}_{0}\coloneqq  v, & \\[2mm]
\ha{v}_{k}\coloneqq  \tbmax{v_{k}}{\J}{\ha{v}_{k-1}} & \text{for $1 \le k \le s$}, \\[2mm]
\kap{\eta}{v}\coloneqq  \ha{v}_{s}, 
\end{cases}
\end{equation}
and then we set
%
%
\begin{equation} \label{eq:zetav}
\zeta(\eta,v)\coloneqq   
\sum_{k=0}^{s-1} \wt (\ha{v}_{k+1} \Rightarrow \ha{v}_{k}). 
\end{equation}

Now, let $i \in I$, and consider the case $\lambda=\vpi_{i}$. 
We fix $N = N_{i} \in \BZ_{\ge 1}$ satisfying the following condition: 
\begin{equation} \label{eq:cN}
N/\pair{\vpi_{i}}{\alpha^{\vee}} \in \BZ \quad 
\text{for all $\alpha \in \Delta^{+}$ such that $\pair{\vpi_{i}}{\alpha^{\vee}} \ne 0$}.
\end{equation}
By the definition of QLS paths of shape $\vpi_{i}$, 
we see that if 
\begin{equation} \label{eq:QLS1}
\eta = (v_{1},\,\dots,\,v_{s} \,;\, 
a_{0},\,a_{1},\,\dots,\,a_{s}) \in \QLS(\vpi_{i}),
\end{equation}
then $Na_{k} \in \BZ$ for all $0 \le k \le s$; 
we write $\eta$ as:
\begin{equation} \label{eq:QLS2}
\eta = ( \underbrace{ v_{1},\dots,v_{1} }_{N(a_{1}-a_{0}) \text{ times}}, 
\underbrace{ v_{2},\dots,v_{2} }_{N(a_{2}-a_{1}) \text{ times}},\,\dots,\,
\underbrace{ v_{s},\dots,v_{s} }_{N(a_{s}-a_{s-1}) \text{ times}}). 
\end{equation}
Let $\eta = (w_{1}, w_{2}, \dots, w_{N}) \in \QLS(\vpi_{i})$ 
(allowing that $w_{k-1}=w_{k}$ for some $2 \le k \le N$), and $v \in W$. 
We define $\ha{w}_{0}\coloneqq  v,\,\ha{w}_{1},\,\ha{w}_{2},\,\dots,\ha{w}_{N}$ by the 
same recursive formula as \eqref{eq:haw}; we see that 
$\ha{w}_{N}=\kappa(\eta,v)$. 
Let $w \in W$. If $\ha{w}_{N} = \kappa(\eta,v) = w$, 
then we deduce from Lemma~\ref{lem:tbmax} that
for each $1 \le k \le N$, there exists a (unique) $\bp_{k} \in \bQBG{\ha{w}_{k}}$ 
such that $\ed(\bp_{k}) = \ha{w}_{k-1}$; note that if $w_{k-1}=w_{k}$, 
then $\ha{w}_{k-1}=\ha{w}_{k}$, and $\bp_{k}$ is 
the trivial directed path $\bt_{\ha{w}_{k}}$. 
%
%
\begin{lem} \label{lem:1/2}
Keep the notation and setting above. 
For $1 \le k \le N$, the directed path $\bp_{k}$ is 
a directed path in $\QBG_{((k-1)/N)\vpi_{i}}(W)$. 
\end{lem}

\begin{proof}
Let $1 \le k \le N$. 
By the definition of QLS paths of shape $\vpi_{i}$, 
there exists a directed path (possibly trivial) from $w_{k}$ to $w_{k-1}$ in 
$\QBG_{((k-1)/N)\vpi_{i}}(\WJ)$. By the same argument as that for \cite[Lemma~6.6\,(2)]{LNSSS2}, 
we deduce that there exists a directed path from $\ha{w}_{k}$ to $\ha{w}_{k-1}$ 
in $\QBG_{((k-1)/N)\vpi_{i}}(W)$. Then, by \cite[Lemma~6.7]{LNSSS2}, 
all shortest directed paths from $\ha{w}_{k}$ to $\ha{w}_{k-1}$ 
in $\QBG(W)$ are directed paths in $\QBG_{((k-1)/N)\vpi_{i}}(W)$. 
In particular, $\bp_{k}$ is a directed path in $\QBG_{((k-1)/N)\vpi_{i}}(W)$.
This proves the lemma. 
\end{proof}

Now, we set 
\begin{equation}
\bQLS{w}: = 
\left\{ \bp=(\bp_{N},\dots,\bp_{2},\bp_{1}) \ \Biggm| \ 
\begin{array}{l}
\text{for all $1 \le k \le N$, $\bp_{k} \in \bQBG{\ed(\bp_{k+1})}$, and}  \\[2mm]
\text{$\bp_{k}$ is a directed path in $\QBG_{((k-1)/N)\vpi_{i}}(W)$}
\end{array} \right\}, 
\end{equation}
where $\bp_{N+1}$ is considered to be the trivial directed path $\bt_{w}$, 
and hence $\ed(\bp_{N+1})=w$. 
Let $\bp=(\bp_{N},\dots,\bp_{2},\bp_{1}) \in \bQLS{w}$, and set 
%
%
\begin{equation} \label{eq:etap}
\eta_{\bp}\coloneqq  (\mcr{\ed(\bp_{2})}^{\J},\dots,\mcr{\ed(\bp_{N})}^{\J},\mcr{\ed(\bp_{N+1})}^{\J}=\mcr{w}^{\J}), 
\end{equation}
where $\J= J_{\varpi_i} = I \setminus \{i\}$. 
By the argument above, we see that 
$\eta_{\bp} \in \QLS(\vpi_{i})$, $\kappa(\eta_{\bp},v) = w$, and $\zeta(\eta_{\bp},v) = 
\sum_{k = 1}^{N} \qwt(\bp_{k})$, with $v = \ed(\bp_{1})$. 
Moreover, we deduce that the map $\bp \mapsto \bigl( \eta_{\bp}, \ed(\bp_{1}) \bigr)$
is a bijection from $\bQLS{w}$ onto 
the set $\bigl\{(\eta,v) \in \QLS(\vpi_{i}) \times W \mid \kappa(\eta,v) =w \bigr\}$.
For $\bp=(\bp_{N},\dots,\bp_{2},\bp_{1}) \in \bQLS{w}$, we set
\begin{equation}
\begin{split}
& \ell(\bp)\coloneqq  \sum_{k = 1}^{N} \ell(\bp_{k}), \qquad
  \ed(\bp)\coloneqq  \ed(\bp_{1}), \\
& \qwt(\bp)\coloneqq  \sum_{k = 1}^{N} \qwt(\bp_{k}), \qquad
  \qwt_{2}(\bp)\coloneqq  \sum_{k = 2}^{N} \qwt(\bp_{k}) = \qwt(\bp)-\qwt(\bp_{1}). 
\end{split}
\end{equation}
%
%
\begin{ex} \label{ex:G2-1}
Assume that $\Fg$ is of type $G_{2}$, and $i = 2$, 
where $\alpha_{2}$ is the long simple root; we have 
$\pair{\alpha_{2}}{\alpha_{1}^{\vee}}=-3$, 
$\pair{\alpha_{1}}{\alpha_{2}^{\vee}}=-1$. 
Note that $\pair{\varpi_2}{\theta^{\vee}} = 2 \not= 1$, 
where $\theta = 3 \alpha_1 + 2 \alpha_2 \in \Delta^{+}$ is the highest root. 
By using Table \eqref{table-G2}, 
the quantum Bruhat graph of type $G_{2}$ is given as in Figure~\ref{fig:G2} below. 
Here, a positive root $\beta \in \Delta^{+}$ is said to be a quantum root if 
$\ell(s_{\beta})=2\pair{\rho}{\beta^{\vee}}-1$ (see also \cite[Section 4.1]{LNSSS1}). 
\begin{equation} \label{table-G2}
\begin{array}{c|c|c|c|c}
\beta \in \Delta^{+} & \text{long or short} & \text{quantum or not} & s_{\beta} 
 & \beta^{\vee} \\ \hline\hline
\alpha_{1} & \text{short} & \text{yes} & s_{1} 
 & \alpha_{1}^{\vee} \\ \hline
\alpha_{1}+\alpha_{2} = s_{2}\alpha_{1} & \text{short} & \text{no} & s_{2}s_{1}s_{2} 
 & \alpha_{1}^{\vee}+3\alpha_{2}^{\vee} \\ \hline
2\alpha_{1}+\alpha_{2} = s_{1}s_{2}\alpha_{1} & \text{short} & \text{no} & s_{1}s_{2}s_{1}s_{2}s_{1} 
 & 2\alpha_{1}^{\vee}+3\alpha_{2}^{\vee} \\ \hline
\alpha_{2} & \text{long} & \text{yes} & s_{2} 
 & \alpha_{2}^{\vee} \\ \hline
3\alpha_{1}+\alpha_{2} = s_{1}\alpha_{2} & \text{long} & \text{yes} & s_{1}s_{2}s_{1} 
 & \alpha_{1}^{\vee}+\alpha_{2}^{\vee} \\ \hline
3\alpha_{1}+2\alpha_{2} = s_{2}s_{1}\alpha_{2} & \text{long} & \text{yes} & s_{2}s_{1}s_{2}s_{1}s_{2} 
 & \alpha_{1}^{\vee}+2\alpha_{2}^{\vee}
\end{array}
\end{equation}

\medskip

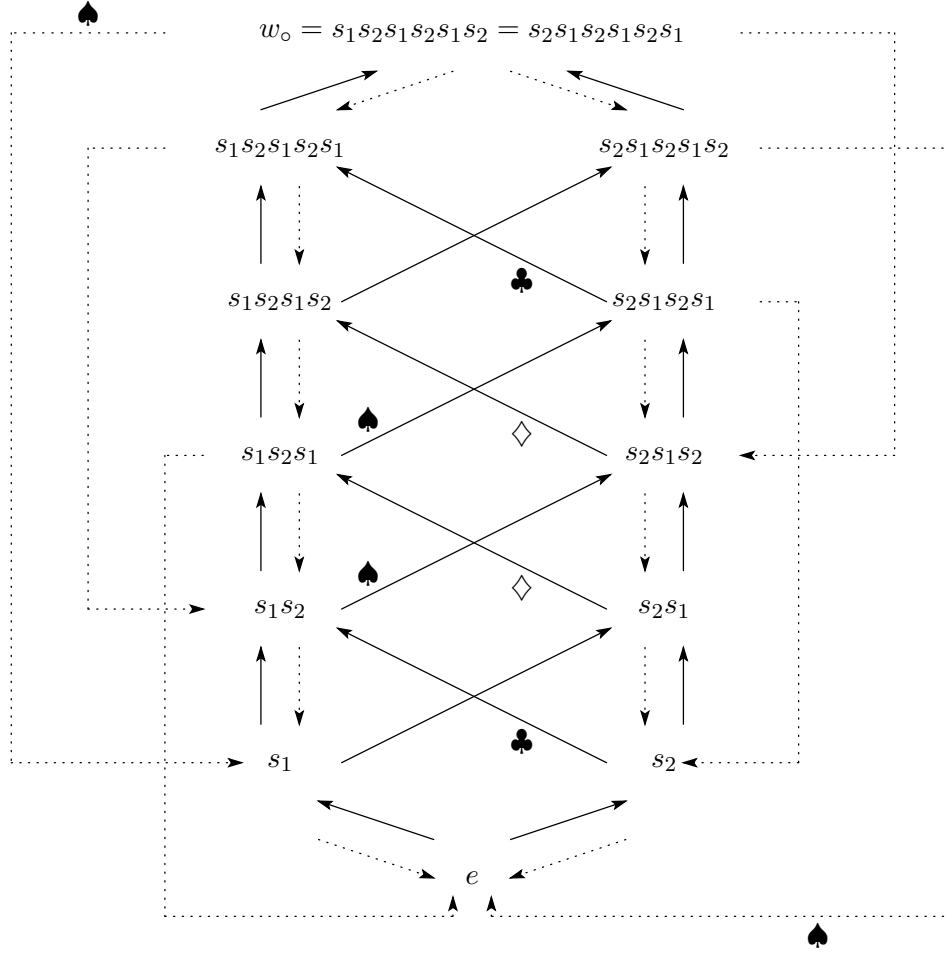
\begin{figure}[ht]
\begin{center}
\input{G2_QBG.tex}
\end{center}

\medskip

\caption{Quantum Bruhat graph of type $G_{2}$. 
We omit the label $\beta$ of an edge $x \edge{\beta} y$ 
if $\pair{\vpi_{i}}{\beta^{\vee}} \in \{0,1\}$. 
The symbol $\spadesuit$ is $3\alpha_{1}+2\alpha_{2}$; 
the edge labeled by $\spadesuit$ is an edge in $\QBG_{(1/2)\vpi_{2}}(W)$. 
The symbols $\diamondsuit$ and $\clubsuit$ are 
$2\alpha_{1}+\alpha_{2}$ and $\alpha_{1}+\alpha_{2}$, respectively; 
the edge labeled by $\diamondsuit$ or $\clubsuit$ is an edge in 
$\QBG_{(1/3)\vpi_{2}}(W) = \QBG_{(2/3)\vpi_{2}}(W)$.} \label{fig:G2}
\end{figure}

In our case, we can take $N=N_{2}=6$. 
We see that 
\begin{align*}
\bp = \bigl( & \bt_{s_{2}}, s_{2} \Be{\alpha_{1}+\alpha_{2}} s_{1}s_{2},
s_{1}s_{2} \Be{3\alpha_{1}+2\alpha_{2}} s_{2}s_{1}s_{2}, \\
& s_{2}s_{1}s_{2} \Be{2\alpha_{1}+\alpha_{2}} s_{1}s_{2}s_{1}s_{2},\bt_{s_{1}s_{2}s_{1}s_{2}},
s_{1}s_{2}s_{1}s_{2} \Be{3\alpha_{1}+\alpha_{2}} s_{2}s_{1}s_{2}s_{1}s_{2} \bigr)
\end{align*}
is an element of $\bQLS{s_{2}}$, and the corresponding QLS path $\eta_{\bp_{1}} \in \QLS(\vpi_{2})$ 
is 
\begin{align*}
\eta_{ \bp } & = (s_{1}s_{2}s_{1}s_{2},s_{1}s_{2}s_{1}s_{2}, s_{2}s_{1}s_{2},s_{1}s_{2}, s_{2},s_{2}) \\
& = \left( s_{1}s_{2}s_{1}s_{2}, s_{2}s_{1}s_{2},s_{1}s_{2}, s_{2} \,;\, 
    0,\frac{1}{3},\frac{1}{2},\frac{2}{3},1 \right). 
\end{align*}
Also, we see that 
\begin{align*}
\bp' & = \bigl( \bt_{w}, \bt_{w}, w \Qe{3\alpha_{1}+2\alpha_{2}} e, \bt_{e}, \bt_{e}, \bt_{e} \bigr), \\
\bp'' & = \bigl( \bt_{w}, \bt_{w}, w \Qe{3\alpha_{1}+2\alpha_{2}} e, \bt_{e}, \bt_{e}, e \Be{\alpha_{2}} s_{2} \bigr)
\end{align*}
are elements of $\bQLS{w}$, with $w \coloneqq   s_{2}s_{1}s_{2}s_{1}s_{2}$, and the corresponding QLS paths are 
\begin{align*}
\eta_{\bp'} = \eta_{\bp''} & = 
(e,e,e,w,w,w) = 
\left( e, w \,;\, 0,\frac{1}{2},1 \right). 
\end{align*}
\end{ex}

%
\subsection{Lakshmibai-Seshadri paths.}
\label{subsec:LS}
%
%
\begin{dfn}[{\cite{Lit94} and \cite{Lit95}}] \label{dfn:LS}
Let $\lambda \in \Lambda^{+}$ be a dominant (integral) weight, and set $\J=\J_{\lambda} = 
\bigl\{ j \in I \mid \pair{\lambda}{\alpha_{j}^{\vee}} = 0 \bigr\}$. 
A  Lakshmibai-Seshadri path (LS path for short) 
of shape $\lambda$ is an element $\eta \in \QLS(\lambda)$ 
of the form \eqref{eq:QLS} satisfying the condition that 
all the edges in a directed path in $\QBG_{a_{k}\nu}(\WJ)$ 
from $w_{k+1}$ to $w_{k}$ are Bruhat edges for all $k = 1,\,2,\,\dots,\,s-1$.
\end{dfn}

Let $\LS(\lambda)$ denote 
the set of all LS paths of shape $\lambda$; 
note that $\LS(\lambda) \subset \QLS(\lambda)$. 
We know from \cite[Lemma~2.14]{MNS} that for $i \in I$, condition
%
%
\begin{equation} \label{eq:QLS=LS}
\QLS(\vpi_{i}) = \LS(\vpi_{i})
\end{equation}
is equivalent to the equality $\pair{\vpi_{i}}{\theta^{\vee}}=1$. 
The following is the list of $i$'s 
satisfying the condition that $\pair{\vpi_i}{\theta^\vee}=1$, 
where our numbering of the nodes of a Dynkin diagram 
is the same as that in \cite[Section 11.4]{H}: 
\begin{center}
\begin{tabular}{l|l}
$A_{n}$ & all $i \in I$ (minuscule) \\ \hline
$B_{n}$ & $i = 1$ (cominuscule), $i = n$ (minuscule), \\
 & with $\alpha_{n}$ the unique short simple root \\ \hline
$C_{n}$ & all $i \in I$ ($\vpi_{1}$ is the unique minuscule weight), \\ 
 & with $\alpha_{n}$ the unique long simple root \\ \hline
$D_{n}$ & $i = 1,\,n-1,\,n$ (minuscule) \\ \hline
$E_{6}$ & $i = 1,\,5$ (minuscule) \\ \hline
$E_{7}$ & $i = 6$ (minuscule)\\ \hline
$E_{8}$ & none \\ \hline
$F_{4}$ & $i = 4$ (not minuscule), with $\alpha_{4}$ a short simple root  \\ \hline
$G_{2}$ & $i = 1$ (not minuscule), with $\alpha_{1}$ the short simple root
\end{tabular}
\end{center}

\begin{rem} \label{rem:012}
Let $i \in I$ be such that $\pair{\vpi_{i}}{\theta^{\vee}}=1$. 
We see that
\begin{equation}
\pair{\vpi_{i}}{\beta^{\vee}} \in \{0,1,2\} \quad 
\text{for all $\beta \in \Delta^{+}$}.
\end{equation}
\end{rem}

%
\subsection{Flag manifolds and Schubert varieties.}

Fix a Borel subgroup $B$ and a parabolic subgroup $P$ of $G$ 
such that $T\subseteq B \subseteq P \subseteq G$.
The opposite Borel subgroup $B^-\subseteq G$ is 
the unique Borel subgroup such that $B\cap B^- = T$. 
The Weyl group $W=\angles*{s_j \mid j \in I}$ of $G$ 
can be identified with $N_G(T)/T$, 
where $N_G(T)$ is the normalizer of $T$ in $G$; 
the Weyl group of $P$, $N_P(T)/T$, is then identified with 
the subgroup $W_{\Kc} \coloneqq  \angles*{s_j \mid j \in \Kc}$ of $W$ 
for a subset $\K \subseteq I$.

Let $Y=G/P$ be the flag manifold associated to $P$. 
Any Weyl group element $w \in W$ defines 
the Schubert varieties $Y_w = \ol{BwP/P}$ and $Y^w=\ol{B^-wP/P}$ in $Y$. 
These varieties depend only on the coset $wW_{\Kc}$ in $W/W_{\Kc}$.
When $w$ belongs to the subset $W^{\Kc} \subseteq W$ of 
 minimal-length coset representatives, we have $\dim Y_w=\codim Y^w=\ell(w)$.

%
\subsection{$K$-theoretic Gromov-Witten invariants on $G/P$.}
\label{subsec:KGW}

For any projective $T$-variety $X$, we denote by $K_T(X)$ 
the Grothendieck group of $T$-equivariant algebraic vector bundles on $X$. 
This ring is an algebra over $K_T(\pt)=R(T)$, 
the representation ring of $T$, which is identified with 
the group algebra $\BZ[\Lambda]$ of $\Lambda$. 
Let $\chi_X: K_T(X)\to K_T(\pt)$ be the pushforward map 
along the structure morphism $X \to \pt$. 
The equivariant $K$-theory ring $K_T(Y)$ of 
the flag manifold $Y=G/P$ has two $K_T(\pt)$-bases 
$\{\cO_w \mid w \in \WK\} = \{\cO_w \mid w \in \WK_{\max}\}$ and $\{\cO^w\mid w \in \WK\}$, 
where $\CO_w=[\CO_{Y_w}]$ and $\CO^w=[\CO_{Y^w}]$ are 
Schubert classes defined by the structure sheaves of 
the Schubert varieties $Y_w$ and $Y^w$, respectively. 
Let $(\CO^w)^{\vee} \in K_T(Y)$ denote the basis of 
$K_T(Y)$ dual to $\CO^w$ for $w \in \WK$, in the sense that 
$\chi^T_Y((\CO^{w})^{\vee} \cdot \CO^{v}) = 
\delta_{w, v}$ for $w, v \in \WK$.
By \cite{Bri}, $(\CO^w)^{\vee}$ is the class of 
the ideal sheaf $\CI_{\partial Y_w}$ of 
the boundary $\partial Y_w$ of $Y_w$.

The homology group $H_{2}(Y; \BZ)$ can be identified 
with $Q^\vee_K$, with $\alpha_j^\vee$ corresponding to 
the class $[Y_{s_j}]$ of the Schubert curve $Y_{s_j}$ for $j\in K$. 
For an effective degree $d \in Q^{\vee,+}_{\K}$ 
and $m\geq 0$, we let $\ol{\CM}_{0,m}(Y,d)$ denote 
the Kontsevich moduli space of $m$-pointed stable maps to 
$Y$ of genus zero and degree $d$ (see \cite{FP}, \cite{Tho}). 
This moduli space is non-empty when $d\neq 0$ or $m\geq 3$. 
In this case, it is equipped with ($T$-equivariant) evaluation maps 
$\ev_k:\ol{\CM}_{0,m}(Y,d)\to Y$ for $1 \leq k\leq m$, 
which send a stable map to the image of the $k$-th marked point in its domain.

For classes 
$\gamma_k \in K_{T}(Y)$, $1 \leq k \leq m$, 
the corresponding $m$-point ($T$-equivariant) 
$K$-theoretic Gromov--Witten (KGW) invariant is defined to be
\begin{equation*}
\langle \gamma_1, \gamma_2, \ldots , \gamma_m \rangle_d^Y \coloneqq   
\chi^{T}_{\ol{\CM}_{0,m}(Y,d)}\left(\prod_{k=1}^m\ev_k^*\gamma_k \right) 
\in K_{T}(\pt).
\end{equation*}
We omit the superscript $Y$ when there is no risk of confusion. 
Non-equivariant KGW invariants are obtained by replacing $T$ 
with the trivial group (which we omit in our notation); 
these KGW invariants take values in $\BZ$. 

%
\subsection{Two-point curve neighborhoods $\Gamma_d(Y_u,Y^{s_i})$.}
\label{sec:Gamma}

Recall the definition of $\Gamma_d(Y_u,Y^{s_i})$ from \eqref{eq:Gamma}; 
if $i \notin K$, i.e., $i \in \Kc$, then we have $s_i \in W_{\Kc}$, and hence $Y^{s_i} = Y^{e} = Y$. 
\begin{lem}\label{lem:Gamma}
Let $d \in Q_{K}^{\vee,+}$, $i \in I$, and $u \in W$. We have
\begin{align*}
  \Gamma_d(Y_u,Y^{s_i}) = 
  \begin{cases}
    \ev_2\parens*{\ev_1^{-1}(Y_u)} & \text{\rm if } d_i>0, \\
    \ev_2\parens*{\ev_1^{-1}(Y_u)}\cap Y^{s_i} & \text{\rm if } d_i=0.
  \end{cases}
\end{align*}
\end{lem}

\begin{proof}
The variety $\Gamma_d(Y_u,Y^{s_i})$ is the closure of 
the union of all rational curves of degree $d$ that intersect both $Y_u$ and $Y^{s_i}$. If $d_i>0$, 
then any such curve intersects $Y^{s_i}$, 
and hence $\Gamma_d(Y_u,Y^{s_i})=\ev_2(\ev_1^{-1}(Y_u))$. 
If $d_i=0$, then any such curve is contained in $Y^{s_i}$, 
and hence $\Gamma_d(Y_u,Y^{s_i})=\ev_2(\ev_1^{-1}(Y_u))\cap Y^{s_i}$. 
This proves the lemma. 
\end{proof}

A geometric argument in \cite[Proposition 3.2(b)]{BCMP1} shows that 
$\ev_2\parens*{\ev_1^{-1}(Y_u)}$ is a $B$-stable Schubert variety; 
see \cite[Theorem 5.1]{BM2}. 
It follows that $\Gamma_d(Y_u,Y^{s_i})$ is a Richardson variety. 
In particular, it is irreducible and has rational singularities. 
Using a result of Koll\'ar \cite{Kol} (see \cite{BCMP1} and \cite{BM1}), 
the aforementioned result \cite[Proposition~3.2]{BCMP1} 
also shows that $\ev_2: \ev_1^{-1}(Y_u) \to \ev_2\parens*{\ev_1^{-1}(Y_u)}$ is 
cohomologically trivial, that is, 
$[\cO_{\ev_2\parens*{\ev_1^{-1}(Y_u)}}]={\ev_2}_*\ev_1^*\cO_u$. 

\begin{cor}
Let $d \in Q_{K}^{\vee,+}$, $i \in I$, and $u \in W$. We have 
\begin{equation}\label{eq:Gamma-cases}
  [\cO_{\Gamma_d(Y_u, Y^{s_i})}] = 
  \begin{cases}
  {\ev_2}_*\ev_1^*\cO_u & \text{\rm if }d_i > 0, \\ 
  ({\ev_2}_*\ev_1^*\cO_u) \cdot \cO^{s_i} & \text{\rm if }d_i=0.
  \end{cases}
\end{equation}
\end{cor}

%
\subsection{Quantum $K$-theory of $G/P$.}

Let $\sQ\coloneqq  \parens*{\sQ_j\mid j\in K}$ be 
the \emph{Novikov} variables. Following \cite{Giv} and \cite{Lee}, 
the $T$-equivariant (small) quantum $K$-theory ring of $Y=G/P$ is  
\begin{equation*} 
QK_T(Y) = K_T(Y) \otimes_{K_T(\pt)} K_T(\pt)\bra{\sQ}
\end{equation*}
as a $K_T(\pt)[\![\sQ]\!]$-module.
It is equipped with a commutative and associative 
product, denoted by $\star$, which is
determined by the condition
\begin{equation} \label{eq:qmulti}
\Kmet{\sigma_1\star\sigma_2}{\sigma_3}=
 \sum_{d \in \QKvp}
 \sQ^d \angles*{\sigma_1,\sigma_2,\sigma_3}_{d} \quad 
 \text{for all  } \sigma_1,\sigma_2,\sigma_3\in K_T(Y),
\end{equation}
where $\sQ^d \coloneqq \prod_{j \in K} \sQ_j^{d_j}$ 
for $d = \sum_{j \in K} d_j \alpha_j^{\vee}$, and 
\begin{equation*}
\Kmet{\sigma_1}{\sigma_2} \coloneqq 
\sum_{d \in \QKvp} \sQ^d \angles*{\sigma_1,\sigma_2}_{d}
\end{equation*}
is the quantum $K$-metric. 

%
\subsubsection{Kato's ring homomorphism.}
\label{subsubsec:kato}

It was proved in \cite{Kat1} and \cite{ACT} that 
for $\sigma_1,\sigma_2\in K_T(Y)$, the product $\sigma_1\star\sigma_2$ 
can always be expressed as a polynomial in $Q$ with coefficients in $K_T(Y)$. 
It follows that 
\begin{equation*}
  QK_T^{\mathrm{poly}}(Y) \coloneqq K_T(Y) \otimes_{K_T(\pt)}K_T(\pt)[\sQ] 
\end{equation*}
is a subring of $QK_T(Y)$. 

Let $\pi: G/B\to G/P$ be the natural map. 
The following theorem is proved in \cite{Kat2} 
(see also \cite[Section~2.3]{KLNS}).

\begin{thm}[Kato]\label{thm:kato-morphism}
There is a surjective ring homomorphism 
\begin{equation*}
    \Psi: QK^\mathrm{poly}_T(G/B)\to QK^\mathrm{poly}_T(G/P),
\end{equation*}
given by $\sigma \mapsto{\pi}_*\sigma$ 
for all $\sigma\in K_T(G/B)$ and $\sQ^d\mapsto\sQ^{\pi_*d}$.
\end{thm}

We shall also write $[d]$ for $\pi_*d$; combinatorially, 
$[\,\cdot\,] \coloneqq [\,\cdot\,]_{\K}$ is 
the projection $Q^\vee \twoheadrightarrow Q^\vee_{\K}$ that 
sends $\sum_{j\in I}c_j\alpha_j^\vee$ to $\sum_{j\in \K}c_j\alpha_j^\vee$.

%
\subsubsection{Chevalley formulas.}
\label{subsec:Monk}

Fix $i \in I$ arbitrarily. 
Let $\lhd$ be a reflection order satisfying \eqref{eq:ro}, 
with $\J=\J_{\vpi_{i}}=I \setminus \{i\}$. 
%
%
\begin{thm}[{\cite[(1.6)]{NOS}; see also \cite[Theorem 49]{LNS}}] \label{thm:NOS}
Let $w \in W$. In $QK_{T}(G/B)$, we have
\begin{equation} \label{eq:NOS}
\CO^{s_{i}} \star \CO^{w} = 
\CO^{w} + 
\sum_{v \in W}
\sum_{
  \begin{subarray}{c}
   \eta \in \QLS(\vpi_{i}) \\
   \kappa(\eta,v) = w
  \end{subarray}}
  (-1)^{\ell(v) - \ell(w) + 1} \be^{-\vpi_{i}+\wt(\eta)}
  \sQ^{\zeta(\eta, v)} \CO^{v}. 
\end{equation}
\end{thm}

By using the bijection $\bQLS{w} \rightarrow 
\bigl\{(\eta,v) \in \QLS(\vpi_{i}) \times W \mid \kappa(\eta,v) =w \bigr\}$, 
$\bp \mapsto \bigl( \eta_{\bp}, \ed(\bp_{1}) \bigr)$, 
given in Section~\ref{subsec:QLS}, we can rewrite equation~\eqref{eq:NOS} 
as follows. 
%
%
\begin{cor} \label{cor:NOS2}
Let $i \in I$ and $w \in W$. In $QK_{T}(G/B)$, we have
\begin{equation} \label{eq:NOS2}
\CO^{s_{i}} \star \CO^{w} = 
\CO^{w} -
\sum_{\bp \in \bQLS{w}} 
  (-1)^{\ell(\bp)}  \be^{-\vpi_{i}+\wt(\eta_{\bp})}
  \sQ^{\qwt(\bp)} \CO^{\ed(\bp)}. 
\end{equation}
\end{cor}
Let $i \in \K$; 
note that $\Kc \subset I \setminus \{i\}=\J$. 
From Theorem~\ref{thm:kato-morphism} and Corollary~\ref{cor:NOS2}, 
we deduce the following Chevalley formula for $QK_T(G/P)$. 

\begin{cor}\label{cor:chev-p}
Let $i \in \K$ and $w \in W^{I\setminus K}$. In $QK_{T}(G/P)$, we have  
\begin{equation} \label{eq:chev-p}
\CO^{s_i} \star \CO^{w} = 
\CO^{w} -
\sum_{\bp \in \bQLS{w}} 
  (-1)^{\ell(\bp)}  \be^{-\vpi_{i}+\wt(\eta_{\bp})}
  {\sQ^{[\qwt(\bp)]}} \CO^{\mcr{\ed(\bp)}}; 
\end{equation}
here, $\mcr{x} = \mcr{x}^{I \setminus K}$ for $x \in W$, 
and $[\,\cdot\,] = [\,\cdot\,]_{\K} : Q^{\vee} \twoheadrightarrow \QKv$ is the projection. 
\end{cor}

\begin{rem}
The formulas in Theorem~\ref{thm:NOS} and Corollary~\ref{cor:NOS2} 
are cancellation-free, but this is not necessarily 
the case for the formula in Corollary~\ref{cor:chev-p}. 
\end{rem}

%
\section{Main Results.}
\label{sec:main}

We follow the notation in Section~\ref{sec:prelim}. 
In particular, let $Y = G/P$ be the flag manifold associated 
to the parabolic subgroup $P \supset B$; the Weyl group of $P$ is 
the subgroup $W_{\Kc} \subset W$ for a subset $K \subset I$. 
Also, $\mcr{x} = \mcr{x}^{\Kc}$, 
$\mxcr{x} = \mxcr{x}^{\Kc}$ for $x \in W$, and 
$[\,\cdot\,] = [\,\cdot\,]_{\K} : Q^{\vee} \twoheadrightarrow \QKv$ is the projection. 
Let $i \in \K$, and $\J = \J_{\vpi_{i}}=I \setminus \{ i \}$.

%
\subsection{Statements of the main results.}
\label{subsec:main}
We now give precise statements of our main results.
Let $d \in Q_{\K}^{\vee,+}$. 
If $i \notin \K$, i.e., $i \in \Kc$, then we have $s_i \in W_{\Kc}$, 
and hence $\CO^{s_i} = [\CO_{Y^{s_i}}] = 1 \in K_T(Y)$; 
in this case, it follows that $\thp{s_i}{w}{x}_{d} = \twp{w}{x}_{d}$
for $w, x \in W$ from the definition of KGW invariants. 
Hence, in the following, we (may and do) assume that $i \in \K$. 
%
%
\begin{thm}\label{thm:parabolic}
Let $i \in \K$, and $w, x \in W$. 
Then, for an effective degree 
$d = \sum_{j \in K} d_j \alpha_j^{\vee} \in \QKvp$ such that $d_i=0$, 
the following holds in $QK_T(Y)${\rm :} 
\begin{equation}\label{eq:parab}
\langle \CO^{s_i}, \CO^{w}, \CO_{x} \rangle_{d} =  
\langle \CO^{s_i} \cdot \CO^{w}, \CO_{x} \rangle_{d},
\end{equation}
where $\CO^{s_i} \cdot \CO^{w}$ denotes the ordinary product in $K_T(Y)$. 
\end{thm}


Let $w \in \WK$, and $x \in \WK_{\max}$. 
For $d=\sum_{j \in \K} d_{j} \alpha_{j}^{\vee} \in \QKvp$, we set 
\begin{equation} \label{eq:pbQLSwxd}
\pbQLS{w,x,d} \coloneqq  
\bigl\{ \bp=(\bp_{N},\dots,\bp_{2},\bp_{1}) \in \bQLS{w} \mid 
[\qwt (\ed(\bp) \Rightarrow x)] \le d - [\qwt(\bp)] \bigr\}; 
\end{equation}
remark that $[\qwt ( \ed(\bp) \Rightarrow x)] \le d - [\qwt(\bp)]$ implies 
$[\qwt(\bp)] \le d$. Also, we define
%
%
\begin{equation} \label{eq:pbR}
\pbR{w,x,d} \coloneqq  
\left\{ \bp=(\bp_{N},\dots,\bp_{2},\bp_{1}) \in \pbQLS{w,x,d} \ \Biggm| \ 
\begin{array}{c}
\pair{\vpi_{i}}{d - [\qwt_{2}(\bp)]} = 0 \\[1.5mm]
\ell(\bp_{1})=0, \, \ed(\bp) \in x\WJs
\end{array}
\right\}.
\end{equation}
When $P=B$, i.e., $\K=I$, we write $\bQLS{w,x,d}$ and $\bR{w,x,d}$ 
for $\pbQLS{w,x,d}$ and $\pbR{w,x,d}$, respectively.
%
%
\begin{thm} \label{thm:qk2p}
Let $i \in K$, and let $w \in \WK$, $x \in \WK_{\max}$. 
Then, for an effective degree $d = \sum_{j \in \K} d_j \alpha_j^{\vee} \in \QKvp$, 
the following holds in $QK_T(Y)${\rm :} 
\begin{equation} \label{eq:qk2pa}
\langle \CO^{s_i}, \CO^w, \CO_x \rangle_{d} =  
 \langle \CO^w, \CO_x \rangle_{d} - 
 \sum_{ \bp \in \pbR{w,x,d} }
(-1)^{\ell(\bp)}\be^{-\vpi_{i}+\wt(\eta_{\bp})}.
\end{equation}
Moreover, the sum on the right-hand side of this equation turns out 
to be zero if $\pair{\vpi_i}{\theta^{\vee}} = 1$ 
(for this condition, see Section~\ref{subsec:LS}) 
and $d_{i} > 0${\rm;} in this case, the following holds in $QK_T(Y)${\rm :} 
\begin{equation} \label{eq:vanishing}
\thp{s_i}{w}{x}_{d} =  \twp{w}{x}_{d}. 
\end{equation}
\end{thm}
%
%
\begin{ex} \label{ex:G2-2}
As in Example~\ref{ex:G2-1}, 
we assume that $\Fg$ is of type $G_{2}$, and $i = 2$; 
note that $\pair{\varpi_{2}}{\theta^{\vee}} \not= 1$, 
where $\theta \in \Delta^{+}$ is the highest root. 
Recall that $\alpha_{2}$ is the long simple root, and that $N=N_{2}=6$. 
Assume that $P = B$, i.e., $\K = I$. 
Let $w \coloneqq s_{2}s_{1}s_{2}s_{1}s_{2}$, and  
$d \coloneqq d_{1}\alpha_{1}^{\vee}+2\alpha_{2}^{\vee} \in Q^{\vee,+}$ 
with $d_{1} \in \BZ_{> 0}$. 
By using \eqref{eq:qk2pa}, we compute $\thp{w}{s_2}{x}_{d}$ for $x \in W$ as follows. 
By Figure~\ref{fig:G2}, we see that 
for all $x \in W$, $\qwt(s_{2}s_{1}s_{2}s_{1}s_{2} \Rightarrow x) \le d$, and hence 
$\twp{s_{2}s_{1}s_{2}s_{1}s_{2}}{x}_{d} = 1$. 
Let $\bp=(\bp_{6},\bp_{5},\dots,\bp_{2},\bp_{1}) \in \bQLS{s_{2}s_{1}s_{2}s_{1}s_{2}}$. 
Then, we deduce that $(\bp_{6},\bp_{5},\dots,\bp_{2})$ is 
\begin{equation*}
\text{either} \quad (\bt_{w},\dots,\bt_{w}) \quad \text{or} \quad 
(\bt_{w},\bt_{w},\bq,\bt_{e},\bt_{e}),
\end{equation*}
where $\bq:w \Qe{3\alpha_{1}+2\alpha_{2}} e$. 
If $\bp \in \bR{w,x,d}$, then 
$(\bp_{6},\bp_{5},\dots,\bp_{2}) = 
(\bt_{w},\bt_{w},\bq,\bt_{e},\bt_{e})$ and $\bp_{1} = \bt_{e}$. 
Note that $\WJs=W_{I \setminus \{2\}} = \{ e,\,s_{1} \}$. Thus, 
we find that 
\begin{equation*}
\bR{w,x,d} = 
 \begin{cases}
 \{ \underbrace{ (\bt_{w},\bt_{w},\bq,\bt_{e},\bt_{e},\bt_{e}) }_{
    \text{$= \bp'$ in Example~\ref{ex:G2-1}} } \} & \text{if $x = e$ or $s_{1}$}, \\[2mm]
 \emptyset & \text{otherwise}. 
 \end{cases}
\end{equation*}
Since $\wt(\eta_{\bp'}) = \vpi_{2}-(3\alpha_{1}+2\alpha_{2}) = 0$, 
we conclude that 
\begin{equation} \label{eq:G2a}
\thp{s_{2}}{w}{x}_{d} =  
 \begin{cases}
 1+\be^{-(3\alpha_1 + 2\alpha_2)} & \text{if $x = e$ or $s_{1}$}, \\
 1 & \text{otherwise}.
 \end{cases}
\end{equation}
\end{ex}

We may also express the correction term in Theorem~\ref{thm:qk2p} 
in terms of QLS paths as follows.
%
%
\begin{prop} \label{prop:qk2p-QLS}
Let $d \in Q_{K}^{\vee,+}$, $i \in K$, and $w, x \in W$ be such that 
$\mcr{w} = w$, $\mxcr{x} = x$. 
We have 
%
%
\begin{equation} \label{eq:pRQLS}
\sum_{ \bp \in \pbR{w,x,d} }
(-1)^{\ell(\bp)}\be^{-\vpi_{i}+\wt(\eta_{\bp})} = 
\sum_{ v \in x\WJs } 
\sum_{
 \begin{subarray}{c}
 \eta \in \QLS(\vpi_{i}) \\[1mm]
 \iota(\eta) \in x\WJs,\,\kappa(\eta,v)=w \\[1mm]
 \pair{\vpi_{i}}{d-[\zeta(\eta,v)]} = 0 \\[1mm]
 [\qwt(v \Rightarrow x)] \le d-[\zeta(\eta,v)]
 \end{subarray} } (-1)^{\ell(v)-\ell(w)} \be^{-\vpi_{i}+\wt(\eta)}. 
\end{equation}
\end{prop}

\begin{proof}
%
Recall from Section~\ref{subsec:QLS} 
the bijection $\bp \mapsto \bigl( \eta_{\bp}, \ed(\bp) \bigr)$ 
from $\bQLS{w}$ onto 
$\bigl\{(\eta,v) \in \QLS(\vpi_{i}) \times W \mid \kappa(\eta,v) =w \bigr\}$; 
notice that $(-1)^{\ell(\bp)} = (-1)^{\ell(\ed(\bp))-\ell(w)}$. 
Hence, it suffices to show that under the bijection above, 
the set $\pbR{w,x,d}$ is mapped to the set
\begin{equation} \label{eq:RQLSa}
\left\{ 
(\eta,v) \in \QLS(\vpi_{i}) \times W \ \Biggm| \ 
\begin{array}{l}
v \in x\WJs,\, \iota(\eta) \in x\WJs,\,\kappa(\eta,v)=w \\[1mm]
\pair{\vpi_{i}}{d-[\zeta(\eta,v)]} = 0 \\[1mm]
[\qwt(v \Rightarrow x)] \le d-[\zeta(\eta,v)]
\end{array}
\right\}. 
\end{equation}
Let $\bp=(\bp_{N},\dots,\bp_{2},\bp_{1}) \in \pbR{w,x,d}$. 
We know that $\kappa(\eta_{\bp},v) = w$ and $\zeta(\eta_{\bp},v) = \qwt(\bp)$, 
with $v = \ed(\bp) \in x\WJs$. Since $\ell(\bp_{1})=0$, it follows that 
$\qwt(\bp) = \qwt_{2}(\bp)$ and $\ed(\bp) = \ed(\bp_{2})$. 
Hence, the equalities $\pair{\vpi_{i}}{d - [\qwt_{2}(\bp)]} = 0$ and 
$[\qwt (\ed(\bp) \Rightarrow x)] \le d - [\qwt(\bp)]$ imply that 
$\pair{\vpi_{i}}{d-[\zeta(\eta_{\bp},v)]} = 0$ and 
$[\qwt(v \Rightarrow x)] \le d-[\zeta(\eta_{\bp},v)]$, respectively. 
Also, we deduce that modulo $\WJs$, $\iota(\eta_{\bp}) \equiv \ed(\bp_{2}) = 
\ed(\bp) = v \equiv x$. 
Thus, $\bigl( \eta_{\bp}, \ed(\bp) \bigr) \in \QLS(\vpi_{i}) \times W$ 
is contained in the set given by \eqref{eq:RQLSa}. 

Let $(\eta,v) \in \QLS(\vpi_{i}) \times W$ be an element 
of the set given by \eqref{eq:RQLSa}, and 
$\bp = (\bp_{N},\dots,\bp_{2},\bp_{1})$ a unique element of 
$\bQBG{w}$ such that $\bigl( \eta_{\bp}, \ed(\bp) \bigr) = (\eta,v)$; 
we will show that $\bp \in \pbR{w,x,d}$. Notice that 
$\ed(\bp) = v \in x\WJs$. Also, since $\ed(\bp_{2}) \equiv 
\iota(\eta_{\bp}) = \iota(\eta) \equiv x$ modulo $\WJs$, it follows that 
$\ed(\bp_{2}) \in x\WJs$. Hence, there exists a shortest directed path 
from $\ed(\bp_{2}) \in x\WJs$ to $\ed(\bp_{1}) \in x\WJs$ in $\QBG(W)$
whose labels are all contained in $\DJs$. Since $\bp_{1}$ is a label-increasing 
(shortest) directed path from $\ed(\bp_{2})$ to $\ed(\bp_{1})$, 
it is lexicographically minimal among all shortest paths from
$\ed(\bp_{2})$ to $\ed(\bp_{1})$; see \cite[Theorem~7.3]{LNSSS1}. 
By \eqref{eq:ro}, we deduce that $\ell(\bp_{1})= 0$. Therefore, we see that 
$[\qwt(\bp)] = [\qwt_{2}(\bp)] = [\zeta(\eta_{\bp},v)] = [\zeta(\eta,v)]$, 
and hence $[\qwt(\ed(\bp) \Rightarrow x)] = [\qwt(v \Rightarrow x)] \le d-[\zeta(\eta,v)]
=d-[\qwt(\bp)]$. Also, we see that $\pair{\vpi_{i}}{d - [\qwt_{2}(\bp)]} = 
\pair{\vpi_{i}}{d - [\zeta(\eta,v)]} = 0$. Thus, we have shown that $\bp \in \pbR{w,x,d}$, as desired. 
This proves the proposition. 
\end{proof}

%
\subsection{Positivity.}
\label{subsec:positivity}

A theorem of Brion~\cite{Bri} states that, 
if a closed irreducible subvariety $Z$ of $Y = G/P$ 
has rational singularities, then the expansion of 
the non-equivariant $K$-theory class $[\cO_Z] \in K(Y)$ 
in the Schubert basis has alternating signs. 
This, combined with Theorems~\ref{thm:parabolic} and \ref{thm:qk2p}, 
has the following consequence.
%
%
\begin{cor}\label{cor:alt-sign}
Let $d \in Q_{K}^{\vee,+}$, $i \in I$, and $u,\,w \in W^{\Kc}$. 
If $d_i=0$ or $\pair{\vpi_i}{\theta^{\vee}} = 1$, 
then the non-equivariant KGW invariants 
$\angles*{\cO^{s_i},\cO_u,\parens*{\cO^w}^\vee}_d$ 
have a positivity property in the sense that 
\begin{equation*}
(-1)^{\ell(w)-\codim\Gamma_d(Y_u,Y^{s_i})}\angles*{\cO^{s_i},\cO_u,\parens*{\cO^w}^\vee}_d\geq 0.
\end{equation*}
\end{cor}

\begin{proof}
Theorems~\ref{thm:parabolic} and \ref{thm:qk2p} imply that 
\begin{align}
\angles*{\cO^{s_i},\cO_u,\parens*{\cO^w}^\vee}_d
= \begin{cases}
  \angles*{\cO_u,\parens*{\cO^w}^\vee}_d 
  & \text{if } d_i>0 \text{ and } \pair{\vpi_i}{\theta^{\vee}} = 1, \\
  \angles*{\cO^{s_i} \cdot \parens*{\cO^w}^\vee,\cO_u}_d 
  & \text{if } d_i=0.
  \end{cases}
\end{align}
Also, we have
\begin{align*}
    \angles*{\cO_u,\parens*{\cO^w}^\vee}_d &= 
    \chi_{\overline{\CM}_{0,2}(Y,d)}\parens*{\ev_1^*\cO_u\cdot\ev_2^*\parens*{\parens*{\cO^w}^\vee}}=
    \chi_Y\parens*{{\ev_2}_*\ev_1^*\cO_u\cdot\parens*{\cO^w}^\vee}\\
& = \chi_Y\parens*{\cO_{\Gamma_d(Y_u,Y^{s_i})}\cdot\parens*{\cO^w}^\vee}\quad\text{if }d_i>0,
\end{align*}
where the second equality follows from the projection formula and 
the third one follows from \eqref{eq:Gamma-cases}, and similarly,
\begin{align*}
    \angles*{\cO^{s_i}\cdot\parens*{\cO^w}^\vee,\cO_u}_d 
& = \chi_{\overline{\CM}_{0,2}(Y,d)}
    \parens*{\ev_2^*\parens*{\cO^{s_i}\cdot{\parens*{\cO^w}^\vee}}\cdot\ev_1^*\cO_u}
  = \chi_Y\parens*{\cO^{s_i}\cdot\parens*{\cO^w}^\vee\cdot{\ev_2}_*\ev_1^*\cO_u}\\ 
& = \chi_Y\parens*{\cO_{\Gamma_d(Y_u,Y^{s_i})}\cdot\parens*{\cO^w}^\vee}\quad\text{if }d_i=0.
\end{align*}
The rest follows from Brion's theorem. This proves the corollary. 
\end{proof}
%
%
\begin{rem}\label{rem:positivity}
Since all relevant maps and classes are $T$-equivariant, 
the above corollary also holds in the equivariant setting 
with the same proof, where we use Anderson, Griffeth, 
and Miller's equivariant generalization of Brion's theorem \cite[Corollary 5.1]{AGM} 
instead. Here, equivariant positivity means 
\begin{equation*}
  (-1)^{\ell(w)-\codim\Gamma_d(Y_u,Y^{s_i})}
  \angles*{\cO^{s_i},\cO_u,\parens*{\cO^w}^\vee}_d\in \BZ_{\geq 0}[\be^{\alpha_j}-1\mid j\in I].
\end{equation*}
\end{rem}

Corollary~\ref{cor:alt-sign} was known for cominuscule varieties $Y = G/P$; 
see \cite[Corollary 4.3]{BCMP2}. The same proof combined 
with geometric statements in \cite{BPX} also imply 
the result for $\SG(2,2n)$, the symplectic Grassmannian of lines. 
In these cases, $\pair{\vpi_i}{\theta^{\vee}} = 1$ is automatic; 
on the other hand, a stronger result was proved, 
where the Schubert divisor $Y^{s_i}$ can be replaced by any Schubert variety $Y^v$.

%
\section{Proofs of Theorems~\ref{thm:parabolic} and \ref{thm:qk2p}.}
\label{sec:prf-main}

%
\subsection{Proof strategy.}
\label{sec:prf-str}

Let us take and fix $i \in \K$. 
We (may and do) assume that $w \in W^{\Kc}$, $x \in W_{\max}^{\Kc}$. 
From equation \eqref{eq:qmulti}, 
we deduce the following equation in $QK_T(Y)$:
\begin{equation}\label{eq:KGWinv}
\sum_{d \in \QKvp} \sQ^d \thp{s_i}{w}{x}_d = 
\Kmet{\CO^{s_i} \star \CO^w}{\CO_x}
\end{equation}
for $w \in \WK$ and $x \in \WK_{\max}$. 
For $v,w\in \WK$, we can always write 
\begin{equation}\label{eq:qmulticoeff}
\CO^{v} \star \CO^w = 
\sum_{\substack{d \in \QKvp \\ z \in W^{I\setminus K}}} N_{v,w}^{z,d} \sQ^d \CO^z
\end{equation}
in $QK_T(Y)$, with $N_{v,w}^{z,d} \in R(T)$.
In the following, we fix  $w \in W^{\Kc}$ arbitrarily, and take 
$\J = J_{\varpi_i} = I \setminus \{i\}$ (see \eqref{eq:J}). 
For $z \in W^{\Kc}$, we set 
\begin{equation}
a_z(\sQ) \coloneqq   \sum_{d \in \QKvp} N_{s_i,w}^{z,d} \sQ^d, 
\end{equation}
which turns out to be an element of $R(T)[\sQ] \subset R(T)\bra{\sQ}$. 
Then, we have in $QK_T(Y)$, 
\begin{equation}
\CO^{s_i} \star \CO^{w} = 
\sum_{z \in W^{\Kc}} a_{z}(\sQ) \CO^{z}; 
\end{equation}
if we set $c_{z}\coloneqq  a_{z}(0) \in R(T)$, then we have in $K_T(Y)$, 
\begin{equation}
\CO^{s_i} \cdot \CO^w = 
\sum_{z \in W^{\Kc}} c_z \CO^z.
\end{equation}
Fix $x \in W^{\Kc}_{\max}$. It follows from equation~\eqref{eq:KGWinv} that 
\begin{equation} \label{eq:main1}
\sum_{d \in \QKvp}
\sQ^{d} \thp{s_{i}}{w}{x}_{d} = 
\sum_{z \in W^{\Kc}} a_{z}(\sQ) \sum_{\xi \in \QKvp} \sQ^{\xi} \twp{z}{x}_{\xi}.  
\end{equation}
When $x\in W^{\Kc}_{\max}$, the 2-point KGW invariants $\twp{z}{x}_{\xi}$ 
can be explicitly described in terms of the quantum Bruhat graph $\QBG(W)$ on $W$ as follows.
%
%
\begin{lem} \label{lem:ptwp}
Let $\xi \in Q^{\vee}_{\K}$, $z \in \WK$, and $x \in \WK_{\max}$. 
Then, we have 
\begin{equation} \label{eq:2ptp-qwt}
\langle \CO^z, \CO_x \rangle_{\xi} = 
\begin{cases}
1 & \text{\rm if $\xi \ge [\qwt(z \Rightarrow x)]$}, \\
0 & \text{\rm otherwise}.
\end{cases}
\end{equation}
\end{lem}

\begin{proof}
It follows from the definition of $\dist_{Y}(z, x)$ in \cite{BCLM} that 
$\dist_{Y}(z, x) = [\dist_X(z, x)]$. Also, from \cite[Lemma~4]{BCLM}, 
we know that $\dist_{X}(z, x) = \qwt(z \Rightarrow x)$ (in our notation). 
Therefore, we deduce that $\dist_{Y}(z, x) = [\qwt(z \Rightarrow x)]$. 
Hence, the assertion of the lemma immediately follows 
from \cite[Proposition~6]{BCLM}. This proves the lemma. 
\end{proof}

We see from Corollary~\ref{cor:chev-p} that for $z \in W^{\Kc}$, 
\begin{equation} \label{eq:N}
a_{z}(\sQ) = \delta_{z,w} - 
  \sum_{ \begin{subarray}{c} \bp \in \bQLS{w} \\ \floor{\ed(\bp)}=z \end{subarray} }(-1)^{\ell(\bp)}
   \be^{-\vpi_{i}+\wt(\eta_{\bp})}  \sQ^{ [\qwt(\bp)] }; 
\end{equation}
since the labels of the edges in $\bp$ are all contained in $\DJp$, we see that
\begin{equation} \label{eq:Qi}
\qwt(\bp) > 0 \iff 
\pair{\vpi_{i}}{\qwt(\bp)} > 0 \iff 
\sQ^{\qwt(\bp)} \in \sQ_{i}R(T)[\sQ]. 
\end{equation}

Our problem is now purely combinatorial. 
We proceed with the proofs of Theorems~\ref{thm:parabolic} and \ref{thm:qk2p} 
in the next two subsections. We also note that Proposition~\ref{prop:comparison} reduces 
the problem from $G/P$ to $G/B$, at least in cases 
where the correction term vanishes; see Appendix~\ref{sec:MX} for more details.

%
\subsection{Proof of Theorem~\ref{thm:parabolic}.}
\label{subsec:proof-parabolic}

Let $i \in \K$, and $w \in W^{\Kc}$, $x \in W_{\max}^{\Kc}$. 
We set 
\begin{equation*}
U \coloneqq  \bigl\{ z \in \WK \mid 
\text{$\mcr{\ed(\bp)} = z$ for some $\bp \in \bQLS{w}$ with $\qwt(\bp) = 0$} \bigr\}. 
\end{equation*}
Notice that $w \in U$ since $(\bt_{w},\dots,\bt_{w},\bt_{w}) \in \bQLS{w}$. 
From \eqref{eq:N} and \eqref{eq:Qi}, we deduce that 
\begin{enu}
\item[\rm (a)] if $z \in U $, then 
$a_{z} ( \sQ ) \in c_{z}  + \sQ_{i}R(T)[ \sQ ]$; 
\item[\rm (b)] if $z \in \WK \setminus U $, then $a_{z} ( \sQ ) \in 
\sQ_{i}R(T)[ \sQ ]$, and hence $c_{z}  = 0$. 
\end{enu}
Indeed, if $z \in U$, then by \eqref{eq:N} and $c_{z} = a_{z}(0)$, we see that 
\begin{equation*}
a_{z}(\sQ) = c_{z} - 
\sum_{ \begin{subarray}{c} \bp \in \bQLS{w} \\ \floor{\ed(\bp)}=z  \\ \qwt(\bp) > 0 \end{subarray} }(-1)^{\ell(\bp)}
   \be^{-\vpi_{i}+\wt(\eta_{\bp})}  \sQ^{ [\qwt(\bp)] }. 
\end{equation*}
From this, it follows that $\sQ^{ [\qwt(\bp)] } \in \sQ_{i}R(T)[\sQ]$ 
for all $\bp \in \bQLS{w}$ such that $\mcr{\ed(\bp)}=z$ and $\qwt(\bp) > 0$ by \eqref{eq:Qi}. 
This proves (a). The proof of (b) is similar; recall that $w \in U$, as seen above. 

Now, let $d = \sum_{j \in \K} d_j \alpha_j^{\vee} \in \QKvp$ be 
such that $d_i = 0$. 
Then, the coefficient of $\sQ^{d}$ on the left-hand side of \eqref{eq:main1} 
is equal to $\thp{s_{i}}{w}{x}_{d}$. 
Also, the right-hand side of \eqref{eq:main1} can be written as: 
\begin{align*}
& 
\sum_{z \in U} \underbrace{a_{z}(\sQ)}_{\in c_{z}+\sQ_{i}R(T)[\sQ]}
\sum_{\xi \in \QKvp} \sQ^{\xi} \twp{z}{x}_{\xi}
+ 
\sum_{z \in W \setminus U}
\underbrace{ a_{z}(\sQ) }_{\in \sQ_{i}R(T)[\sQ] } 
\sum_{\xi \in \QKvp} \sQ^{\xi} \twp{z}{x}_{\xi} \\[3mm]
& = 
\sum_{z \in U} c_{z} 
\sum_{ \begin{subarray}{c}
 \xi = \sum_{j \in I} \xi_{j}\alpha_{j}^{\vee} \in \QKvp \\
 \xi_{i} = 0
 \end{subarray} } \sQ^{\xi} \twp{z}{x}_{\xi} + \text{(an element in $\sQ_{i}R(T)[\sQ]$)}.
\end{align*}
Since $d_{i}= 0$ (and hence $\sQ^{d} \notin \sQ_{i}R(T)[\sQ]$) by the assumption, 
we see that the coefficient of $\sQ^{d}$ 
on the right-hand side of \eqref{eq:main1} is equal to 
$\sum_{z \in U} c_{z} \twp{z}{x}_{d}$.
Therefore, we conclude that 
\begin{equation*}
\thp{s_{i}}{w}{x}_{d} = \sum_{z \in U} c_{z} \twp{z}{x}_{d} 
\stackrel{\text{(b)}}{=}
\sum_{z \in W^{\Kc}} c_{z} \twp{z}{x}_{d}, 
\end{equation*}
as desired. This proves Theorem~\ref{thm:parabolic}. 

%
\subsection{Proof of Theorem~\ref{thm:qk2p}.}
\label{subsec:proof-pcorrection}

Let $i \in K$, and $w \in W^{\Kc}$, $x \in W_{\max}^{\Kc}$. 
By \eqref{eq:N}, we can rewrite \eqref{eq:main1} as follows:
\begin{align} 
& \sum_{d \in Q_{K}^{\vee,+}}
  \sQ^{d} \thp{s_{i}}{w}{x}_{d} \nonumber \\
& = \sum_{\xi \in Q_{K}^{\vee,+}} \sum_{z \in \WK}
\left( \delta_{z,w} \sQ^{\xi}\twp{z}{x}_{\xi} - 
  \sum_{ \begin{subarray}{c} \bp \in \bQLS{w} \\ \mcr{\ed(\bp)}=z \end{subarray} }
  (-1)^{ \ell(\bp) } \be^{-\vpi_{i}+\wt(\eta_{\bp})} \sQ^{ [\qwt(\bp)] + \xi }\twp{z}{x}_{\xi}\right) \nonumber \\
& = \sum_{\xi \in Q_{K}^{\vee,+}} \sum_{z \in \WK} \delta_{z,w} \sQ^{\xi}\twp{z}{x}_{\xi} \nonumber \\
& \hspace*{20mm} - 
\sum_{\xi \in Q_{K}^{\vee,+}} 
\sum_{ \bp \in \bQLS{w} } (-1)^{\ell(\bp)} 
\be^{-\vpi_{i}+\wt(\eta_{\bp})}
\sQ^{ [\qwt(\bp)] + \xi }\twp{\mcr{\ed(\bp)}}{x}_{\xi}. \label{eq:main2}
\end{align}
By comparing the coefficients of $\sQ^{d}$ on the leftmost-hand side and the rightmost-hand side, 
we deduce that 
%
%
\begin{align}
\thp{s_{i}}{w}{x}_{d}
& = \twp{w}{x}_{d} - 
\sum_{ \begin{subarray}{c} \bp \in \bQLS{w} \\ [\qwt(\bp)] \le d \end{subarray} } 
(-1)^{\ell(\bp)}
\be^{-\vpi_{i}+\wt(\eta_{\bp})}
\twp{\mcr{\ed(\bp)}}{x}_{d-[\qwt(\bp)]} \nonumber \\
& = \twp{w}{x}_{d} - 
\sum_{ \begin{subarray}{c} \bp \in \bQLS{w} \\ 
[ \qwt(\ed(\bp) \Rightarrow x)]  \le d-[\qwt(\bp)] \end{subarray} }
(-1)^{\ell(\bp)}\be^{-\vpi_{i}+\wt(\eta_{\bp})}, \nonumber
\end{align}
where for the second equality, we use Lemma~\ref{lem:ptwp}. 
Since $[ \qwt(\mcr{\ed(\bp)} \Rightarrow x) ] = [ \qwt(\ed(\bp) \Rightarrow x) ]$ 
by \cite[Lemma~7.2]{LNSSS2}, we find that 
%
%
\begin{equation} \label{eq:qk2p-2}
\langle \CO^{s_i} , \CO^w , \CO_x  \rangle_{d} 
= \langle \CO^w , \CO_x  \rangle_{d} - 
\sum_{ \bp \in \pbQLS{w,x,d} }
(-1)^{\ell(\bp)}\be^{-\vpi_{i}+\wt(\eta_{\bp})}.
\end{equation}
%
%
We set 
\begin{align}
\pbQLS{w,x,d,+} & = \bigl\{ \bp \in \pbQLS{w,x,d} \mid \pair{\vpi_{i}}{d-\qwt_{2}(\bp)} > 0 \bigr\}, \\
\pbQLS{w,x,d,0} & = \bigl\{ \bp \in \pbQLS{w,x,d} \mid \pair{\vpi_{i}}{d-\qwt_{2}(\bp)} = 0 \bigr\}; 
\end{align}
note that $\pbQLS{w,x,d} = \pbQLS{w,x,d,+} \sqcup \pbQLS{w,x,d,0}$, and 
$\pbR{w,x,d} \subset \pbQLS{w,x,d,0}$. 
In order to prove Theorem~\ref{thm:qk2p}, 
it suffices to prove the following proposition and 
that $ \pbR{w,x,d}=\emptyset$ when $\pair{\vpi_i}{\theta^\vee}=1$ and $d_i>0$.
\begin{prop}\label{prop:qk2p}
  We have 
\begin{equation} \label{eq:qk2-a}
 \sum_{ \bp \in \pbQLS{w,x,d,+} }
(-1)^{\ell(\bp)}\be^{-\vpi_{i}+\wt(\eta_{\bp})} = 0, 
\end{equation}
\begin{equation} \label{eq:qk2-b}
 \sum_{ \bp \in \pbQLS{w,x,d,0} \setminus \pbR{w,x,d} }
(-1)^{\ell(\bp)}\be^{-\vpi_{i}+\wt(\eta_{\bp})} = 0. 
\end{equation}
\end{prop}

To prove Proposition~\ref{prop:qk2p}, we will first consider 
the case $\K=I$ for the simpler notation, 
and then point out how to modify the proof for general $\K \subset I$.

\subsubsection{Proof of \eqref{eq:qk2-a}.}

First, assume that $P = B$, i.e., $K=I$. 
In this case, $W^{\Kc} = W^{\Kc}_{\max} = W$, $\QKv = Q^{\vee}$; 
$\mcr{x} = \mxcr{x} = x$ for all $x \in W$, and 
the projection $Q^{\vee} \twoheadrightarrow \QKv$ is just the identity map. 

In the proof below, we give two ``sijections'' (i.e., bijections 
between signed sets) $\Theta$ and $\Theta'$.

First, for $\bp=(\bp_{N},\dots,\bp_{2},\bp_{1}) \in \bQLS{w}$, 
we define $\Theta(\bp)$ as follows. Here we recall that 
$\FL{\bp_{1}}$ denotes the final label of $\bp_{1}$ 
(if $\bp_{1}$ is not the trivial directed path); 
if $\bp_{1}$ is the trivial directed path, then we set $\FL{\bp_{1}}\coloneqq  -\infty$, 
which is a formal element such that 
$- \infty \lhd \alpha$ for all $\alpha \in \Delta^{+}$. 
Remark that $\alpha_{i}$ is the maximum element in $\Delta^{+}$ with respect to $\lhd$. 
\begin{enu}
\item If $\FL{\bp_{1}}=\alpha_{i}$, then 
we define $\Theta(\bp_{1})$ to be the directed path obtained 
from $\bp_{1}$ by removing the final edge $\ed(\bp_{1})s_{i} \edge{\alpha_{i}} \ed(\bp_{1})$; 
we set $\Theta(\bp)\coloneqq  (\bp_{N},\dots,\bp_{2},\Theta(\bp_{1}))$. 
\item If $\FL{\bp_{1}} \ne \alpha_{i}$, then 
we define $\Theta(\bp_{1})$ to be the directed path obtained from $\bp_{1}$ 
by adding the edge $\ed(\bp_{1}) \edge{\alpha_{i}} \ed(\bp_{1})s_{i}$ at the end of $\bp_{1}$; 
we set $\Theta(\bp)\coloneqq  (\bp_{N},\dots,\bp_{2},\Theta(\bp_{1}))$. 
\end{enu}
It is easily verified that 
$\Theta(\bp) \in \bQLS{w}$, with $\ell(\Theta(\bp))=\ell(\bp) \pm 1$, and that 
$\eta_{\Theta(\bp)} = \eta_{\bp}$. Hence, $\Theta$ is a sijection on $\bQLS{w}$; 
however, the set $\bQLS{w,x,d}$ is not stable under $\Theta$ in general. 

\medskip

Now, we set 
\begin{align*}
\bA & \coloneqq   \bigl\{ (\bp_{N},\dots,\bp_{2},\bp_{1}) \in \bQLS{w,x,d,+} \mid 
  \FL{\bp_{1}} = \alpha_{i} \bigr\}, \\
\bB & \coloneqq   \bigl\{ (\bp_{N},\dots,\bp_{2},\bp_{1}) \in \bQLS{w,x,d,+} \mid 
  \FL{\bp_{1}} \ne \alpha_{i} \bigr\}. 
\end{align*}
%
%
\begin{lem} \label{lem:Theta1}
We have $\Theta(\bp) \in \bB$ for all $\bp \in \bA$.
\end{lem}

\begin{proof}
Let $\bp = (\bp_{N},\dots,\bp_{2},\bp_{1}) \in \bA$, and set $z\coloneqq  \ed(\bp)$; 
recall that $\qwt(z \Rightarrow x) \le d - \qwt(\bp)$ and 
$\pair{\vpi_{i}}{d-\qwt_{2}(\bp)} > 0$. 
We see that $\qwt_{2}(\Theta(\bp)) = \qwt_{2}(\bp)$, which implies that 
$\pair{\vpi_{i}}{d-\qwt_{2}(\Theta(\bp))} = \pair{\vpi_{i}}{d-\qwt_{2}(\bp)} > 0$. 
Note that 
$\ed(\Theta(\bp)) = zs_{i} \edge{\alpha_{i}} z=\ed(\bp)$ is the final edge of $\bp_{1}$. 
Let $zs_{i} \edge{\alpha_{i}} z \Rightarrow x$ be 
the concatenation of the edge $zs_{i} \edge{\alpha_{i}} z$ 
with a shortest directed path from $z$ to $w$. 
Assume that $zs_{i} \edge{\alpha_{i}} z$ is a Bruhat edge. 
We see that $\qwt(\Theta(\bp))=\qwt(\bp)$, and that 
\begin{align*}
\qwt(\ed(\Theta(\bp)) \Rightarrow x) & = 
\qwt(zs_{i} \Rightarrow x) \le 
\qwt(zs_{i} \Be{\alpha_{i}} z \Rightarrow x) = \qwt(z \Rightarrow x) \\ 
& \le d-\qwt(\bp) = d-\qwt(\Theta(\bp)); 
\end{align*}
for the first inequality, see Proposition~\ref{prop:qwt}. 
Thus, we deduce that 
$\qwt(\ed(\Theta(\bp)) \Rightarrow x) \le d-\qwt(\Theta(\bp))$, 
and hence $\Theta(\bp) \in \bB$. 
Assume that $zs_{i} \edge{\alpha_{i}} z$ is a quantum edge. 
We see that $\qwt(\Theta(\bp))=\qwt(\bp)-\alpha_{i}^{\vee}$, and that 
\begin{align*}
\qwt(\ed(\Theta(\bp)) \Rightarrow x) & = 
\qwt(zs_{i} \Rightarrow x) \le 
\qwt(zs_{i} \Qe{\alpha_{i}} z \Rightarrow x) = \alpha_{i}^{\vee}+ \qwt(z \Rightarrow x) \\ 
& \le \alpha_{i}^{\vee} + d-\qwt(\bp) = d-\qwt(\Theta(\bp)). 
\end{align*}
Thus, we deduce that 
$\qwt(\ed(\Theta(\bp)) \Rightarrow x) \le d-\qwt(\Theta(\bp))$, 
and hence $\Theta(\bp) \in \bB$. This proves the lemma. 
\end{proof}

We divide the subset $\bB$ into two subsets: 
\begin{align*}
\bB_{1} & \coloneqq   \bigl\{ \bp \in \bB \mid \Theta(\bp) \in \bQLS{w,x,d,+} \bigr\}, \\
\bB_{2} & \coloneqq   \bigl\{ \bp \in \bB \mid \Theta(\bp) \not\in \bQLS{w,x,d,+} \bigr\}.
\end{align*}
By Lemma~\ref{lem:Theta1}, we find that $\bp \mapsto \Theta(\bp)$ 
is a sijection on the set $\bA \sqcup \bB_{1}$. 
We need to define another sijection $\Theta'$ on $\bB_{2}$. 
Let $\bp = (\bp_{N},\dots,\bp_{2},\bp_{1}) \in \bB_{2}$, and set $z \coloneqq   \ed(\bp)$. 
Since $\qwt_{2}(\Theta(\bp)) = \qwt_{2}(\bp)$, it follows that 
$\pair{\vpi_{i}}{d-\qwt_{2}(\Theta(\bp))} = \pair{\vpi_{i}}{d-\qwt_{2}(\bp)} > 0$. 
Hence, we have
\begin{equation} \label{eq:notle}
\qwt ( \ed(\Theta(\bp)) \Rightarrow x) \not\le d - \qwt(\Theta(\bp)).
\end{equation}
Let $\sT_{1} \in \{\sB, \sq\}$ be the type of the edge 
$\ed(\bp)=z \edge{\alpha_{i}} zs_{i}$ (the final edge of $\Theta(\bp_{1})$), and set 
\begin{equation} \label{eq:osT}
\ol{\sT}_{1}\coloneqq  \begin{cases} 
 \sB & \text{if $\sT_{1}=\sq$}, \\
 \sq & \text{if $\sT_{1}=\sB$}.
 \end{cases}
\end{equation}
Note that $\qwt(\Theta(\bp)) = \qwt(\bp) + \delta_{\sT_{1},\bq}\alpha_{i}^{\vee}$. 
Since $zs_{i} \te{\alpha_{i}}{\ol{\sT}_{1}} z \Rightarrow x$ 
(the concatenation of the edge $zs_{i} \edge{\alpha_{i}} z$ 
 with a shortest directed path from $z$ to $w$) 
is a directed path from $zs_{i}$ to $x$, we have 
\begin{align*}
& \qwt(\ed(\Theta(\bp)) \Rightarrow x) = 
\qwt(zs_{i} \Rightarrow x) \\
& = 
\begin{cases}
-\delta_{\sT_1,\sq} \alpha_{i}^{\vee} + 
\qwt(z \Rightarrow x) & \text{if condition (D) is satisfied}, \\
\delta_{\sT_1,\sB} \alpha_{i}^{\vee} + \qwt(z \Rightarrow x) 
& \text{if condition (D) is not satisfied},
\end{cases}
\end{align*}
where condition (D) is: 
\begin{enu}
\item[(D)] the initial label $\IL{\dec(z \Rightarrow x)}$ of 
the label-decreasing directed path $\dec(z \Rightarrow x)$ 
from $z$ to $x$ is  $\alpha_{i}$; 
recall that $\alpha_{i}$ is the maximum element with respect to $\lhd$. 
\end{enu}
Recall that $\qwt(\ed(\bp) \Rightarrow x) \le d - \qwt(\bp)$. 
If condition (D) were satisfied, then we would have 
$\qwt(\ed(\Theta(\bp)) \Rightarrow x) \le d - \qwt(\Theta(\bp))$, 
which contradicts \eqref{eq:notle}. 
Therefore, condition (D) is not satisfied in this case.
Now, we set 
\begin{equation} \label{eq:beta}
\beta\coloneqq   \begin{cases}
 \FL{\bp_{1}} & \text{if $\bp_{1}$ is not the trivial one}, \\[1.5mm]
 -\infty & \text{if $\bp_{1}$ is the trivial one}; 
 \end{cases}
\end{equation}
note that 
$\beta \ne \alpha_{i}$ since $\bp \in \bB_{2} \subset \bB$. 
Also, we set 
\begin{equation} \label{eq:gamma}
\gamma\coloneqq   \begin{cases}
 \IL{\dec(z \Rightarrow x)} & \text{if $z \ne x$}, \\[1.5mm]
 -\infty & \text{if $z = x$}; 
 \end{cases}
\end{equation}
note that $\gamma \ne \alpha_{i}$ since condition (D) is not satisfied. 
%
%
\begin{lem} \label{lem:s1}
Keep the notation and setting above. 
If $\gamma \in \DJs \cup \{ -\infty \}$, then $\qwt(\bp_{1}) > 0$, 
and hence $\bp_{1}$ is not the trivial directed path. 
\end{lem}
%
%
\begin{rem} \label{rem:bg}
By Lemma~\ref{lem:s1}, we see that $(\beta,\gamma) \ne (-\infty,-\infty)$. 
Recall that if there exists a directed path of the form 
$z's_{\beta'} \edge{\beta'} z' \edge{\beta'} z's_{\beta'}$ in $\QBG(W)$ 
for some $z' \in W$ and $\beta' \in \Delta^{+}$, then 
$\beta'$ is a simple root. 
Since $\alpha_{i}$ is the unique simple root contained in 
$\DJp$, and since $\beta,\,\gamma \ne \alpha_{i}$ as seen above, 
we conclude that $\beta \ne \gamma$. 
\end{rem}

\begin{proof}[Proof of Lemma~\ref{lem:s1}]
Suppose, for a contradiction, that 
$\gamma \in \DJs \cup \{-\infty\}$ and $\qwt(\bp_{1}) = 0$; 
note that $\qwt(\bp)=\qwt_{2}(\bp)$ in this case. 
Recall that $\qwt(\Theta(\bp)) = \qwt(\bp) + \delta_{\sT_{1},\sq}\alpha_{i}^{\vee}$. 
Also, we deduce by \eqref{eq:ro} that 
$\qwt(z \Rightarrow x) \in \sum_{j \in \J} \BZ_{\ge 0}\alpha_{j}^{\vee}$. 
Since $\pair{\vpi_{i}}{d-\qwt_{2}(\bp)} > 0$, and 
$\qwt(z \Rightarrow x) \le d - \qwt(\bp) = d - \qwt_{2}(\bp)$, we see that 
\begin{align*}
\qwt(\ed(\Theta(\bp)) \Rightarrow x) & = 
\delta_{\sT_{1},\sB} \alpha_{i}^{\vee} + \qwt(z \Rightarrow x) \\
& \le d - \qwt(\bp) - \delta_{\sT_{1},\sq}\alpha_{i}^{\vee} = d-\qwt(\Theta(\bp)), 
\end{align*}
which contradicts \eqref{eq:notle}. 
This proves the lemma. 
\end{proof}

We define $\Theta'(\bp_{1})$ as follows (recall that $\beta \ne \gamma$, and 
$\beta,\,\gamma \ne \alpha_{i}$). 
\begin{enu}
\item[(i)] If $\beta \rhd \gamma$, then 
$\Theta'(\bp_{1})$ is defined to be the directed path obtained from $\bp_{1}$ by 
removing the final edge (labeled by $\beta$). 

\item[(ii)] If $\beta \lhd \gamma$, then 
$\Theta'(\bp_{1})$ is defined to be the directed path obtained from $\bp_{1}$ by 
adding the edge labeled by $\gamma$ at the end of $\bp_{1}$ 
(note that if $\beta=-\infty$, then $\gamma \in \DJp$; see Lemma~\ref{lem:s1}). 
\end{enu}
We set $\bq_{1}\coloneqq  \Theta'(\bp_{1})$. We see that 
$\Theta'(\bp)\coloneqq   (\bp_{N},\dots,\bp_{2},\Theta'(\bp_{1})) \in \bQLS{w}$; notice that 
\begin{equation} \label{eq:sij}
\ell(\Theta'(\bp)) = \ell(\bp) \pm 1 \quad \text{and} \quad 
\eta_{\Theta'(\bp)}=\eta_{\bp}.
\end{equation}
We claim that $\bq \coloneqq   \Theta'(\bp) \in \bB_{2}$. 
We remark that the final label $\FL{\bq_{1}}$ of 
$\bq_{1}=\Theta'(\bp_{1})$ is not equal to $\alpha_{i}$ (recall the definition of $\bB$). 
Since $\qwt_{2}(\Theta'(\bp)) = \qwt_{2}(\bp)$, it follows that 
$\pair{\vpi_{i}}{d-\qwt_{2}(\Theta'(\bp))} = \pair{\vpi_{i}}{d-\qwt_{2}(\bp)} > 0$. 
Hence, it remains to show that
\begin{equation} \label{eq:B1}
 \underbrace{ \qwt(\ed(\bq) \Rightarrow x) }_{= \, \qwt(\ed(\bq_{1}) \Rightarrow x)} 
 \le d-\qwt(\bq), 
\end{equation}
\begin{equation} \label{eq:B2}
 \underbrace{ \qwt( \ed(\Theta(\bq) ) \Rightarrow x) }_{= \, \qwt(\ed(\Theta(\bq_{1})) \Rightarrow x)} 
 \not \le d-\qwt(\Theta(\bq)). 
\end{equation}

\medskip

\paragraph{\bf Case 1.}
%
Assume that $\beta \rhd \gamma$ (see (i));
note that $\ed(\bq) = \ed(\bq_{1}) = zs_{\beta} \edge{\beta} z=\ed(\bp_{1}) = \ed(\bp)$ 
is the final edge of $\bp_{1}$. Let $\sT_{\beta} \in \{\sB, \sq\}$ 
be the type of the edge $zs_{\beta} \edge{\beta} z$. We have 
\begin{equation*}
\underbrace{
\overbrace{\ed(\bp_{2}) \edge{\bullet} \cdots \edge{\bullet}
\ed(\bq_{1})}^{= \, \bq_{1} \, = \, \Theta'(\bp_{1})} = zs_{\beta} \te{\beta}{\sT_{\beta}} 
\ed(\bp_{1})}_{= \, \bp_{1}} = 
\underbrace{z \edge{\gamma} \cdots \edge{\bullet} x}_{= \, \dec(z \Rightarrow x)}, 
\end{equation*}
\begin{equation*}
\qwt(\Theta'(\bp)) = \qwt(\bp) - \delta_{\sT_{\beta},\sq} \beta^{\vee}, \qquad 
\qwt(z \Rightarrow x) = \qwt(zs_{\beta} \Rightarrow x) - \delta_{\sT_{\beta},\sq} \beta^{\vee}. 
\end{equation*}
Therefore, it follows that 
\begin{align*}
\qwt(\ed(\Theta'(\bp)) \Rightarrow x) & = 
  \qwt(\ed(\bp) \Rightarrow x) + \delta_{\sT_{\beta},\sq} \beta^{\vee} \\
& \le d-\qwt(\bp) + \delta_{\sT_{\beta},\sq} \beta^{\vee} = d - \qwt(\Theta'(\bp)), 
\end{align*}
which implies \eqref{eq:B1}. Let us show \eqref{eq:B2}. 
Let $\sT_{2} \in \{\sB, \sq\}$ be the type of the edge 
$zs_{\beta} \edge{\alpha_{i}} zs_{\beta}s_{i}$, and set
\begin{equation} \label{eq:osT2}
\ol{\sT}_{2}\coloneqq  \begin{cases} 
 \sB & \text{if $\sT_{2} = \sq$}, \\
 \sq & \text{if $\sT_{2} = \sB$}.
 \end{cases}
\end{equation}
Notice that
\begin{equation*}
zs_{\beta}s_{i} \te{\alpha_{i}}{\ol{\sT}_{2}} zs_{\beta} \te{\beta}{\sT_{\beta}} 
\underbrace{z \edge{\gamma} \cdots \edge{\bullet} x}_{ = \, \dec(z \Rightarrow x) }, \qquad
zs_{i} \te{\alpha_{i}}{\ol{\sT}_{1}}
\underbrace{z \edge{\gamma} \cdots \edge{\bullet} x}_{ = \, \dec(z \Rightarrow x) }
\end{equation*}
are both label-decreasing directed paths (recall that $\alpha_{i} \rhd \beta \rhd \gamma$). 
Hence, we see that
\begin{equation*}
\begin{split}
\qwt(\Theta(\bq)) 
 & = \qwt(\Theta(\bp)) -\delta_{\sT_{1},\sq} \alpha_{i}^{\vee} 
     -\delta_{\sT_{\beta},\sq} \beta^{\vee} + \delta_{\sT_{2},\sq}\alpha_{i}^{\vee}, \\
\qwt(zs_{\beta}s_{i} \Rightarrow x) & = \qwt(zs_{i} \Rightarrow x) 
  - \delta_{\ol{\sT}_{1},\bq} \alpha_{i}^{\vee} 
  + \delta_{\sT_{\beta},\sq} \beta^{\vee} 
  + \delta_{\ol{\sT}_{2},\sq} \alpha_{i}^{\vee} \\
  & = \qwt(zs_{i} \Rightarrow x) 
  - (-\alpha_{i}^{\vee}+\delta_{\ol{\sT}_{1},\bq} \alpha_{i}^{\vee})
  + \delta_{\sT_{\beta},\sq} \beta^{\vee} 
  + (-\alpha_{i}^{\vee}+\delta_{\ol{\sT}_{2},\sq} \alpha_{i}^{\vee}) \\
  & = \qwt(zs_{i} \Rightarrow x) 
  + \delta_{\sT_{1},\bq} \alpha_{i}^{\vee}
  + \delta_{\sT_{\beta},\sq} \beta^{\vee} 
  - \delta_{\sT_{2},\sq} \alpha_{i}^{\vee}.
\end{split}
\end{equation*}
Since $\qwt(zs_{i} \Rightarrow x) \not\le d - \qwt(\Theta(\bp))$ by the assumption, 
it follows that $\qwt(zs_{\beta}s_{i} \Rightarrow x) \not\le d - \qwt(\Theta(\bq))$. 
Thus, we have shown \eqref{eq:B2}, as desired. 

\medskip

\paragraph{\bf Case 2.}
%
Assume that $\beta \lhd \gamma$ (see (ii));
note that $\ed(\bp) = \ed(\bp_{1}) = z \edge{\gamma} zs_{\gamma} = \ed(\bq_{1}) = \ed(\bq)$ 
is the initial edge of $\dec(z \Rightarrow x)$. Let $\sT_{\gamma} \in \{\sB, \sq\}$ 
be the type of the edge $z \edge{\gamma} zs_{\gamma}$. We have 
\begin{equation*}
\underbrace{
\ed(\bp_{2}) \edge{\bullet} \cdots \edge{\bullet} zs_{\beta} \edge{\beta} 
\ed(\bp_{1})}_{= \, \bp_{1}} = 
\underbrace{
z \te{\gamma}{\sT_{\gamma}} \ed(\bq_{1}) = zs_{\gamma} \edge{\bullet} \cdots \edge{\bullet} x}_{%
= \, \dec(z \Rightarrow x)}, 
\end{equation*}
\begin{equation*}
\qwt(\Theta'(\bp)) = \qwt(\bp) + \delta_{\sT_{\gamma},\sq}\gamma^{\vee}, \qquad 
\qwt(z \Rightarrow x) = \qwt(zs_{\gamma} \Rightarrow x) + \delta_{\sT_{\gamma},\sq}\gamma^{\vee}. 
\end{equation*}
Therefore, it follows that 
\begin{align*}
\qwt(\ed(\Theta'(\bp)) \Rightarrow x) 
& = \qwt(\ed(\bp) \Rightarrow x) -  \delta_{\sT_{\gamma},\sq}\gamma^{\vee} \\
& \le d-\qwt(\bp) - \delta_{\sT_{\gamma},\sq}\gamma^{\vee} = d-\qwt(\Theta'(\bp)), 
\end{align*}
which implies \eqref{eq:B1}. Let us show \eqref{eq:B2}.
Let $\sT_{3} \in \{\sB, \sq\}$ be the type of the edge 
$zs_{\gamma} \edge{\alpha_{i}} zs_{\gamma}s_{i}$, and 
set
\begin{equation} \label{eq:osT3}
\ol{\sT}_{3}\coloneqq  \begin{cases} 
 \sB & \text{if $\sT_{3} = \sq$}, \\
 \sq & \text{if $\sT_{3} = \sB$}.
 \end{cases}
\end{equation}
Notice that 
\begin{equation*}
zs_{i} \te{\alpha_{i}}{\ol{\sT}_{1}}
\underbrace{z \te{\gamma}{\sT_{\gamma}} zs_{\gamma} \edge{\bullet} 
\cdots \edge{\bullet} x}_{ = \, \dec(z \Rightarrow x) }, \qquad
zs_{\gamma}s_{i} \te{\alpha_{i}}{\ol{\sT}_{3}}
\underbrace{zs_{\gamma} \edge{\bullet} \cdots \edge{\bullet} x}_{ \text{label-decreasing} }
\end{equation*}
are both label-decreasing directed paths. Hence, we see that
\begin{align*}
\qwt(\Theta(\bq)) & = \qwt(\Theta(\bp)) 
 - \delta_{\sT_{1},\sq}\alpha_{i}^{\vee}
 + \delta_{\sT_{\gamma},\sq}\gamma^{\vee} 
 + \delta_{\sT_{3},\sq}\alpha_{i}^{\vee}, \\
\qwt(zs_{\gamma}s_{i} \Rightarrow x) 
& = \qwt(zs_{i} \Rightarrow x) 
    - \delta_{\ol{\sT}_{1},\sq} \alpha_{i}^{\vee} 
    - \delta_{\sT_{\gamma},\sq}\gamma^{\vee} 
    + \delta_{\ol{\sT}_{3},\sq}\alpha_{i}^{\vee} \\
& = \qwt(zs_{i} \Rightarrow x) 
    - (-\alpha_{i}^{\vee} + \delta_{\ol{\sT}_{1},\sq} \alpha_{i}^{\vee})
    - \delta_{\sT_{\gamma},\sq}\gamma^{\vee} 
    + (-\alpha_{i}^{\vee} + \delta_{\ol{\sT}_{3},\sq}\alpha_{i}^{\vee}) \\
& = \qwt(zs_{i} \Rightarrow x) 
    + \delta_{\sT_{1},\sq}\alpha_{i}^{\vee}
    - \delta_{\sT_{\gamma},\sq}\gamma^{\vee} 
    - \delta_{\sT_{3},\sq}\alpha_{i}^{\vee}. 
\end{align*}
Since $\qwt(zs_{i} \Rightarrow x) \not\le d - \qwt(\Theta(\bp))$ by the assumption, 
it follows that $\qwt(zs_{\gamma}s_{i} \Rightarrow x) \not\le d - \qwt(\Theta(\bq))$. 
Thus, we have shown \eqref{eq:B2}, as desired. 

\medskip

We can easily deduce that $\Theta'(\Theta'(\bp)) = \bp$ for all $\bp \in \bB_{2}$. 
By this fact, together with \eqref{eq:sij}, $\Theta'$ is a sijection on $\bB_{2}$. 
Therefore, we conclude that 
\begin{align*}
& \sum_{\bp \in \bQLS{w,x,d,+}} (-1)^{\ell(\bp)} \be^{-\vpi_{i}+\wt(\eta_{\bp})} \\[2mm]
& = 
\underbrace{ \sum_{\bp \in \bA \sqcup \bB_{1}} (-1)^{\ell(\bp)} \be^{-\vpi_{i}+\wt(\eta_{\bp})} }_{= \, 0 \ 
\text{by sijection $\Theta$}}  + 
\underbrace{ \sum_{\bp \in \bB_{2}} (-1)^{\ell(\bp)} \be^{-\vpi_{i}+\wt(\eta_{\bp})} }_{= \, 0 \
\text{by sijection $\Theta'$}} = 0,
\end{align*}
as desired.

\medskip

For general $K$, we set 
\begin{align*}
\bA^{\K} & := \bigl\{ (\bp_{N},\dots,\bp_{2},\bp_{1}) \in \pbQLS{w,x,d,+} \mid 
  \FL{\bp_{1}} = \alpha_{i} \bigr\}, \\
\bB^{\K} & := \bigl\{ (\bp_{N},\dots,\bp_{2},\bp_{1}) \in \pbQLS{w,x,d,+} \mid 
  \FL{\bp_{1}} \ne \alpha_{i} \bigr\}, \\
\bB^{\K}_{1} & := \bigl\{ \bp \in \bB^{\K} \mid \Theta(\bp) \in \pbQLS{w,x,d,+} \bigr\}, \\
\bB^{\K}_{2} & := \bigl\{ \bp \in \bB^{\K} \mid \Theta(\bp) \not\in \pbQLS{w,x,d,+} \bigr\}. 
\end{align*}
By the same arguments as that for Lemma~\ref{lem:Theta1} (and that after its proof), 
we can verify that $\Theta(\bp) \in \bB^{\K}_{1}$ for all $\bp \in \bA^{\K}$, 
and hence that $\Theta$ is a sijection on the set $\bA^{\K} \sqcup \bB^{\K}_{1}$. 
Also, a sijection $\Theta'$ on $\bB^{\K}_{2}$ can be defined in exactly the same way as the one 
on $\bB_{2}$ as follows. Let $\bp = (\bp_{N},\dots,\bp_{2},\bp_{1}) \in \bB^{\K}_{2}$, 
and set $z := \ed(\bp)$. We define $\beta$ and $\gamma$ in the same manner as 
in \eqref{eq:beta} and \eqref{eq:gamma}, respectively.
As in Lemma~\ref{lem:s1} and Remark~\ref{rem:bg}, we deduce that 
$\beta \ne \alpha_{i}$, $\gamma \ne \alpha_{i}$, and $\beta \ne \gamma$. 
Then, we define $\Theta'(\bp_{1})$ in exactly the same way as in (i) and (ii), 
and set $\Theta'(\bp):= (\bp_{N},\dots,\bp_{2},\Theta'(\bp_{1}))$. 
By replacing $\qwt(\bullet)$, $\beta^{\vee}$, $\gamma^{\vee}$ in Cases 1 and 2 by $[\qwt(\bullet)]$, $[\beta^{\vee}]$, $[\gamma^{\vee}]$, 
respectively, we can show that $\Theta'(\bp) \in \bB^{\K}_{2}$. 
Thus, we have defined a sijection $\Theta'$ on $\bB^{\K}_{2}$. 
By using the sijections $\Theta$ (on $\bA^{\K} \sqcup \bB^{\K}_{1}$) and 
$\Theta'$ (on $\bB^{\K}_{2}$), we can prove equation \eqref{eq:qk2-a}, as desired. 

\subsubsection{Proof of \eqref{eq:qk2-b}.}
%
First, assume that $P = B$, i.e., $K=I$.
In this case, $W^{\Kc} = W^{\Kc}_{\max} = W$, $\QKv = Q^{\vee}$; 
$\mcr{x} = \mxcr{x} = x$ for all $x \in W$, and the projection $Q^{\vee} \twoheadrightarrow \QKv$ is just the identity map. 

For simplicity of notation, we set
\begin{equation}
\bS{w,x,d}\coloneqq  \bQLS{w,x,d,0} \setminus \bR{w,x,d}. 
\end{equation}
Let $\bp = (\bp_{N},\dots,\bp_{2},\bp_{1}) \in \bS{w,x,d}$, and set $z\coloneqq  \ed(\bp)$. 
We define $\beta$ and $\gamma$ as in \eqref{eq:beta} and \eqref{eq:gamma}, respectively. 
Since $\bp \in \bS{w,x,d}$, we see that
\begin{equation} \label{eq:bg2}
\text{if $\gamma \in \DJs \cup \{-\infty\}$, then $\beta \in \DJp$}. 
\end{equation}
Suppose, for a contradiction, that $(\beta,\gamma) = (\alpha_{i},\alpha_{i})$. 
Then we have $zs_{i} \edge{\beta=\alpha_{i}} z \edge{\gamma=\alpha_{i}} zs_{i}$; 
note that either of these two edges is a quantum edge. 
Recall that $\qwt_{2}(\bp)+\qwt(\bp_{1}) = \qwt(\bp) \le d$. 
Since $\pair{\vpi_{i}}{d-\qwt_{2}(\bp)} = 0$, and 
since the labels of the edges in $\bp_{1}$ are all contained in $\DJp$, 
we deduce that $\qwt(\bp_{1})=0$, or equivalently, all the edges in $\bp_{1}$ 
are Bruhat edges; in particular, $zs_{i} \edge{\beta=\alpha_{i}} z$ is a Bruhat edge. 
Also, we have $\qwt(\bp) = \qwt_{2}(\bp)$. 
Since $\qwt(z \Rightarrow x) \le d - \qwt(\bp) = d-\qwt(\bp_{2})$, 
it follows that an edge in $\dec(z \Rightarrow x)$ whose label 
is contained in $\DJp$ is a Bruhat edge; in particular, 
$z \edge{\gamma=\alpha_{i}} zs_{i}$ is also a Bruhat edge. However, 
this is a contradiction. Thus, we have shown that $(\beta,\gamma) \ne (\alpha_{i},\alpha_{i})$. 
By the same argument as that in Remark~\ref{rem:bg}, we deduce that $\beta \ne \gamma$. 
Now, we define $\Psi(\bp_{1})$ as follows. 
\begin{enu}
\item[(iii)] If $\beta \rhd \gamma$, then 
$\Psi(\bp_{1})$ is defined to be the directed path obtained from $\bp_{1}$ by 
removing the final edge (labeled by $\beta$). 

\item[(iv)] If $\beta \lhd \gamma$, then 
$\Psi(\bp_{1})$ is defined to be the directed path obtained from $\bp_{1}$ by 
adding the edge labeled by $\gamma$ at the end of $\bp_{1}$ 
(note that if $\beta=-\infty$, then $\gamma \in \DJp$; see \eqref{eq:bg2}). 
\end{enu}
Then we set $\Psi(\bp)\coloneqq  (\bp_{N},\dots,\bp_{2},\Psi(\bp_{1}))$; 
notice that $\Psi(\bp) \in \bQLS{w}$, and that 
\begin{equation}
\ell(\Psi(\bp)) = \ell(\bp) \pm 1 \quad \text{and} \quad 
\eta_{\Psi(\bp)} = \eta_{\bp}.
\end{equation}
We claim that $\Psi(\bp) \in \bS{w,x,d}$. 
Since $\qwt_{2}(\Psi(\bp)) = \qwt_{2}(\bp)$, 
it follows that $\pair{\vpi_{i}}{d-\qwt_{2}(\Psi(\bp))}=
\pair{\vpi_{i}}{d-\qwt_{2}(\bp)}=0$. We see that 
$\ell(\Psi(\bp_1)) > 0$ or $\ed(\Psi(\bp)) \notin x\WJs$. 
Hence, it remains to show that 
\begin{equation} \label{eq:Psi}
\qwt(\ed(\Psi(\bp)) \Rightarrow x) \le d - \qwt(\Psi(\bp)).
\end{equation}
By the argument in the previous paragraph, we deduce that 
in both cases (iii) and (iv), 
$\qwt(\ed(\Psi(\bp)) \Rightarrow x) = \qwt(\ed(\bp) \Rightarrow x)$ and 
$\qwt(\Psi(\bp)) = \qwt(\bp)$. Since 
$\qwt(\ed(\bp) \Rightarrow x) \le d - \qwt(\bp)$ by the assumption, 
we obtain \eqref{eq:Psi}, as desired. 

Notice that $\ell(\Psi(\bp)) = \ell(\bp) \pm 1$ and 
$\eta_{\Psi(\bp)}=\eta_{\bp}$. Hence, $\Psi$ is a sijection on $\bS{w,x,d}$. 
Thus, we have proved \eqref{eq:qk2-b} in the $K=I$ case.

\medskip

For general $K$, we set
\begin{equation}
\pbS{w,x,d}:=\pbQLS{w,x,d,0} \setminus \pbR{w,x,d}. 
\end{equation}
Let $\bp = (\bp_{N},\dots,\bp_{2},\bp_{1}) \in \pbS{w,x,d}$, and set $z:=\ed(\bp)$. 
Let $\beta$ and $\gamma$ be as above. Then, we see that \eqref{eq:bg2} in this case, too. 
Also, by the same argument as that after \eqref{eq:bg2}, 
we deduce that $(\beta,\gamma) \ne (\alpha_{i},\alpha_{i})$, and hence $\beta \ne \gamma$. 
We now define $\Psi(\bp_{1})$ in exactly the same way as in (iii) and (iv), 
and set $\Psi(\bp):= (\bp_{N},\dots,\bp_{2},\Psi(\bp_{1}))$. 
By replacing $\qwt(\bullet)$ in the proof above 
by $[\qwt(\bullet)]$, we can show that $\Psi(\bp) \in \pbS{w,x,d}$. 
Thus, we have defined a sijection $\Psi$ on $\pbS{w,x,d}$. 
By using the sijection $\Psi$, we can prove \eqref{eq:qk2-b}, as desired. 

\subsubsection{Vanishing of correction terms.} 
\label{subsubsec:vanish}

Here we assume that $\pair{\vpi_{i}}{\theta^{\vee}}=1$ and $d_{i} > 0$. 
By \eqref{eq:QLS=LS}, we have $\QLS(\vpi_{i})=\LS(\vpi_{i})$. 
Also, by Remark~\ref{rem:012}, we can take $N=N_{i}=2$. Note that 
if $\eta = (w_{1},w_{2}) \in \QLS(\vpi_{i})=\LS(\vpi_{i})$ 
(with the notation as in \eqref{eq:QLS2}), then 
$w_{1} \ge w_{2}$ in the Bruhat order $\ge$. 
%

We will show that the set $\pbR{w,x,d}$, given by \eqref{eq:pbR}, 
is empty; this implies that 
$\thp{s_i}{w}{x}_{d} =  
 \twp{w}{x}_{d}$ in $QK_T(Y)$.
Indeed, we show that $\pair{\vpi_{i}}{d-[\qwt_{2}(\bp)]} > 0$
for all $\bp=(\bp_{2},\bp_{1}) \in \pbQLS{w,x,d}$. 
Recall that $\eta_{\bp}=(\mcr{\ed(\bp_{2})}^{\J}, \mcr{w}^{\J}) \in \QLS(\vpi_{i}) 
= \LS(\vpi_{i})$. 
As seen above, we have $\mcr{\ed(\bp_{2})}^{\J} \ge \mcr{w}^{\J}$ 
in the Bruhat order $\ge$ on $\WJ$, where $\J=\J_{\vpi_{i}}=I \setminus \{i\}$.
Here, by Lemma~\ref{lem:tbmax}, we have $w = \tbmax{w}{\J}{\ed(\bp_{2})}$. 
Therefore, by Lemma~\ref{lem:tbmax-B}, the inequality 
$\mcr{\ed(\bp_{2})}^{\J} \ge \mcr{w}^{\J}$ implies that 
$\ed(\bp_{2}) \ge w$ in the Bruhat order. 
Hence, we deduce that $\qwt_{2}(\bp) = \qwt(\bp_{2}) = 0$. 
Thus, we conclude that $\pair{\vpi_{i}}{d-[\qwt_{2}(\bp)]} = 
\pair{\vpi_{i}}{d} = d_{i} > 0$, as desired. 
This completes the proof of Theorem~\ref{thm:qk2p}. 

\appendix

%
\section{A Peterson comparison formula in quantum $K$-theory \newline (Notes by Mihalcea and Xu).}
\label{sec:MX}

In Appendix~\ref{sec:MX}, we state a comparison formula for 
the $K$-theoretic GW invariants on $G/P$'s, generalizing the similar formula 
in cohomology, conjectured by Peterson \cite{Pet97} and proved by Woodward \cite{Woo}. 
Most of our arguments can be found, implicitly or explicitly, in Woodward's proof, but we indicate here the 
modifications needed to work in $K$-theory. In particular, the key fact that for $m\geq 1$, the natural map 
$\Mb_{0,m}(G/B, \dl{d}) \to \Mb_{0,m}(G/P, d)$ induced by the natural projection of flag manifolds $G/B \to G/P$ is 
cohomologically trivial, may be deduced directly from arguments in Woodward's proof. 
We use the comparison formula to reduce Theorems~\ref{thm:parabolic} and \ref{thm:qk2p} from $G/P$ to $G/B$. 

Recall that locally, the irreducible projective variety $\Mb_{0,m}(G/P, d)$ may be realized as a smooth space modulo 
a finite group (see, e.g., \cite{FP}), and therefore it has rational singularities \cite{Bou}. 
Let $\Hom_{d}(G/P)$ denote the set of morphisms $\BP^1 \to G/P$ of degree $d$. 
It is known that when non-empty, it is a dense smooth open subset 
of the moduli space $\Mb_{0,3}(G/P, d)$ of $3$-point, degree $d$ maps to $G/P$; 
see \cite{Tho,FP}. 

For $K\subseteq J$ subsets of $I$, let $P$ and $Q$ be parabolic subgroups of $G$ containing $B$ 
corresponding to the subsets $I\setminus \K$ and $I \setminus \J$, respectively. 
We shall write $\pi_{\J \setminus \K} : G/Q \to G/P$ for the natural projection. 

We know from \cite[Lemma~1]{Woo} 
(see also \cite[Remark~10.17]{LS} and \cite[Lemma~3.8]{LNSSS1}) 
that for each $d \in Q^{\vee}_K \cong H_2(G/P, \BZ)$, 
there exists a unique $\dl{d} \in Q^{\vee} \cong H_2(G/B, \BZ)$ 
which is sent to $d$ under the map $H_2(G/B, \BZ) \to H_2(G/P, \BZ)$ induced by the natural projection $G/B \to G/P$ 
(combinatorially, the map $H_2(G/B, \BZ) \to H_2(G/P, \BZ)$ corresponds to the projection $[\,\cdot\,] = [\,\cdot\,]_\K : 
Q^{\vee} \twoheadrightarrow Q^{\vee}/Q^{\vee}_{\Kc} \cong Q^{\vee}_{K}$, 
given by $[\sum_{j \in I} d_j \alpha_j^{\vee}] = \sum_{j \in K} d_j \alpha_j^{\vee}$) and 
such that $\pair{\alpha}{\dl{d}} \in \{ 0, -1 \}$ holds 
for all positive roots in $\Delta^+ \cap Q_{\Kc}$, 
with $Q_{\Kc} \coloneqq   \sum_{j \in \Kc} \BZ \alpha_j$; 
the element $\dl{d}$ is called the \emph{Peterson lift} of $d$. 
It was stated by Peterson and proved by Woodward in \cite[Lemma 1]{Woo} 
that if $\Hom_{d}(G/P)$ is non-empty, so is $\Hom_{\dl{d}}(G/B)$. 
This $\dl{d}$ determines a parabolic subgroup $P^{\prime} \subseteq P$ containing $B$ corresponding to the subset
$I \setminus K^{\prime} \coloneqq   \{j\in \Kc \mid \pair{\alpha_j}{\dl{d}}=0\}$.
Let $d^{\prime} \coloneqq   \left(\pi_{I\setminus K^{\prime}}\right)_*\dl{d} \in H_2(G/P^{\prime}; \BZ)$.

The main result proved by Woodward \cite[Theorem 3]{Woo} is the following.

\begin{thm}[{\cite{Woo}}] \label{thm:woodward}
The following hold. 
\begin{enumerate}
\item The morphism $\mathrm{Hom}_{\dl{d}}(G/B) \to \mathrm{Hom}_{d^{\prime}}(G/P^{\prime}) \times_{G/P^{\prime}} G/B$ given by $f \mapsto (\pi_{I\setminus K^{\prime}} \circ f, f(0))$ is an open, dense immersion{\rm ;} 
\item The morphism $\mathrm{Hom}_{d^{\prime}}(G/P^{\prime}) \to \mathrm{Hom}_{d}(G/P)$ given by $f\mapsto\pi_{K^{\prime} \setminus K} \circ f$ is birational. 
\end{enumerate}
\end{thm} 

A morphism $f: X \to Y$ of irreducible varieties is said to be \emph{cohomologically trivial} if $f_* \cO_X = \cO_Y$ and $R^i f_* \cO_X = 0$ for $i>0$. A result of Koll\'ar \cite{Kol} shows that if $f$ is surjective, $X$ and $Y$ are projective with rational singularities, and the general fibers of $f$ are rationally connected, then $f$ is cohomologically trivial; we refer to \cite[Theorem 3.1]{BM1} and \cite[Proposition 5.2]{BCMP1} for details.  This, together with the previous theorem, implies the following. 

\begin{cor}\label{cor:cohtriv} For $m\geq 0$, the natural morphisms 
\[ \Mb_{0,m+1}(G/B,\dl{d}) \to \Mb_{0,m+1}(G/P^{\prime}, d^{\prime}) \times_{G/P^{\prime}} G/B, \quad 
   \Mb_{0,m}(G/P^{\prime},d^{\prime}) \to \Mb_{0,m}(G/P, d)\]
are birational, surjective, and cohomologically trivial.
\end{cor}
\begin{proof}
The $m\geq 2$ case for the first morphism and the $m\geq 3$ case 
for the second morphism follow from Theorem~\ref{thm:woodward}.

For small $m$, consider the following commutative diagrams of surjective morphisms:
\begin{equation*}
\begin{CD}
\ol{\CM}_{0,3}(G/B, \ha{d}) @>>> \ol{\CM}_{0,3}(G/P',d')\times_{G/P'}G/B \\
@VVV @VVV \\
\ol{\CM}_{0,m+1}(G/B, \ha{d}) @>>> \ol{\CM}_{0,m+1}(G/P',d')\times_{G/P'}G/B, 
\end{CD}
\end{equation*}
\begin{equation*}
\begin{CD}
\ol{\CM}_{0,3}(G/P', d') @>>> \ol{\CM}_{0,3}(G/P,d) \\
@VVV @VVV \\
\ol{\CM}_{0,m}(G/P', {d}') @>>> \ol{\CM}_{0,m}(G/P,d).
\end{CD}
\end{equation*}
In each case, the general fiber of the morphism on the bottom row has dimension $0$ 
and is connected because the general fibers of the other three morphisms are connected.
Taking also into account that all the varieties considered have rational singularities,  
the results follow. 
\end{proof}

\begin{rem} This corollary also appears on \cite[page 8]{Woo}, and it leads to the main calculation of the Peterson comparison formula.\end{rem}
\begin{cor}\label{cor:cohtriv-final} 
For $m\geq 0$, the natural map $\Mb_{0,m}(G/B, \dl{d}) \to \Mb_{0,m}(G/P, d)$ is cohomologically trivial.\end{cor}
\begin{proof} This follows from Corollary~\ref{cor:cohtriv} and \cite[Lemma~2.4]{BCMP2}.
{Indeed, we see that the morphism $\Mb_{0,m+1}(G/B, \dl{d}) \to \Mb_{0,m}(G/P, d)$ is cohomologically trivial by using the sequence of maps 
\begin{equation*}
\Mb_{0,m+1}(G/B, \dl{d})
\to \Mb_{0,m+1}(G/P', d') \to \Mb_{0,m+1}(G/P, d) \to \Mb_{0,m}(G/P, d), 
\end{equation*}
noting that the forgetful
maps are cohomologically trivial. Then the claim follows by using
the sequence of maps
$\Mb_{0,m+1}(G/B, \dl{d}) \to \Mb_{0,m}(G/B, \dl{d}) \to \Mb_{0,m}(G/P, d)$.}
This proves the corollary. 
\end{proof}

Let $\pi:X=G/B \to Y=G/P$ denote the natural projection, and 
$\Phi: \Mb_{0,m}(X, \dl{d}) \to \Mb_{0,m}(Y, d)$ the induced map. 
\begin{prop}\label{prop:Kcomparison} For $a_1,\dots,a_m \in K_T(Y)$, and $d \in H_2(Y; \mathbb{Z})$ an effective degree, there is an equality of KGW invariants{\rm :} 
\[ \langle a_1,\dots,a_m \rangle^{Y}_{d} = \langle \pi^*a_1 ,\dots, \pi^*a_m \rangle_{\dl{d}}^X \/. \]
\end{prop}
\begin{proof} We have a commutative diagram: 
\begin{equation*}
\begin{CD}
\Mb_{0,m}(X,\dl{d}) @>{EV}>> X^m \\
@V{\Phi}VV @ VV{\pi^m}V \\
\Mb_{0,m}(Y,d) @>{EV_{Y}}>> Y^m, 
\end{CD}
\end{equation*}
where $EV$ and $EV_Y$ are the products of the evaluation maps into $X$ and $Y$, respectively, at each of the $m$ marked points. 
Starting from the definition of KGW invariants, we see that 
\[ \begin{split} \langle a_1,\dots,a_m \rangle^{Y}_{d} & = \chi^T_{\Mb_{0,m}(Y,d)}\left(EV_Y^*\prod_{k=1}^m a_k \cdot [\cO_{\Mb_{0,m}(Y,d)}]\right) \\ 
& = \chi^T_{\Mb_{0,m}(Y,d)}\left(EV_Y^*\prod_{k=1}^m a_k  \cdot \Phi_*[\cO_{\Mb_{0,m}(X,\dl{d})}]\right) \\
& = \chi^T_{\Mb_{0,m}(X,\dl{d})}\left(\Phi^* EV_Y^*\prod_{k=1}^m a_k  \cdot [\cO_{\Mb_{0,m}(X,\dl{d})}]\right) \\
& = \chi^T_{\Mb_{0,m}(X,\dl{d})}\left(EV^*\prod_{k=1}^m\pi^*a_k \cdot [\cO_{\Mb_{0,m}(X,\dl{d})}]\right) \\
& = \langle \pi^*a_1 ,\dots, \pi^*a_m \rangle_{\dl{d}}^X \/. 
\end{split}
\]
Here, the second equality follows from the cohomological triviality in Corollary~\ref{cor:cohtriv-final}, {the third one follows from the projection formula \cite[Chap.~5, Sect.~5.3.12]{CG}, and the fourth one is just using that the diagram commutes.} 
This proves the proposition. 
\end{proof} 

\begin{cor}
If equation \eqref{eq:parab} in Theorem~\ref{thm:parabolic} (resp., equation \eqref{eq:vanishing} in Theorem~\ref{thm:qk2p}) holds for $X$, then it also holds for $Y$.
\end{cor} 

\begin{proof}
  Note that for $i \in K$, $d_i=0$ if and only if $\dl{d}_i=0$. 
  
  Using Proposition~\ref{prop:Kcomparison}, Theorems \ref{thm:parabolic} and \ref{thm:qk2p} for $X$, we see that 
  \begin{align*}
    &\langle \CO^{s_i} , \CO^u,  \CO_v  \rangle_{d}^{Y}=\angles*{\pi^*\CO^{s_i} ,\pi^*\CO^u,\pi^* \CO_v }^X_{\dl{d}}=\angles*{\CO^{s_i},\pi^*\CO^u,\pi^* \CO_v }_{\dl{d}}^X,
  \end{align*}
  and when $d_i>0$ and $\pair{\vpi_i}{\theta^{\vee}} = 1$,
  \begin{align*}
    \angles*{\CO^{s_i},\pi^*\CO^u,\pi^* \CO_v }_{\dl{d}}^X=\angles*{\pi^*\CO^u,\pi^* \CO_v }_{\dl{d}}^X=\langle \CO^{u}, \CO_{v} \rangle_{d}^Y; 
  \end{align*}
  when $d_i=0$, we see that 
  \begin{align*}
    &\angles*{\CO^{s_i},\pi^*\CO^u,\pi^* \CO_v }_{\dl{d}}^X = \angles*{\pi^*\CO^{s_i} \cdot\pi^*\CO^u,\pi^* \CO_v }_{\dl{d}}^X = \angles*{\pi^*\CO^{s_i} \cdot\pi^*\CO^u,\pi^* \CO_v ,\cO_X}_{\dl{d}}^X \\ &= \angles*{\pi^*(\CO^{s_i} \cdot\CO^u),\pi^* \CO_v ,\pi^*\cO_Y}_{\dl{d}}^X = \angles*{\CO^{s_i} \cdot\CO^u, \CO_v ,\cO_Y}_{d}^Y = \angles*{\CO^{s_i} \cdot\CO^u, \CO_v }_{d}^Y,
  \end{align*}
  where we have used the fact that $\pi^*\cO_Y=\cO_X$ for the natural projection $\pi:X\to Y$. This proves equations \eqref{eq:parab} and \eqref{eq:vanishing} for $Y$. 
\end{proof}

Also, using similar arguments, we can deduce the following from Theorem~\ref{thm:qk2p}. 
\begin{prop}
  \[
    \sum_{ \bp \in \bR{w,x,\dl{d}} }
    (-1)^{\ell(\bp)}\be^{-\vpi_{i}+\wt(\eta_{\bp})}=\sum_{ \bp \in \pbR{w,x,d} }
    (-1)^{\ell(\bp)}\be^{-\vpi_{i}+\wt(\eta_{\bp})}. 
  \]
\end{prop}

It would be interesting to find a purely combinatorial proof of this proposition.

%
\section{Relationship between bases.}
\label{sec:basis}

In Appendix~\ref{sec:basis}, we assume that $P = B$, i.e., $K = I$; 
recall that $X = G/B$. 
Let $(\CO^x)^{\vee} \in K_T(X)$, $x \in W$, 
denote the basis of $K_T(X)$ dual to $\CO^x$, $x \in W$, 
in the sense that $\chi^T_X((\CO^{x})^{\vee} \cdot \CO^{y}) = 
\delta_{x, y}$ for $x, y \in W$, 
where $\chi^T_X$ is the ($T$-equivariant) pushforward 
along the structure morphism of $X$.
In fact, for $x \in W$, $(\CO^x)^{\vee}$ is the class of 
the ideal sheaf $\CI_{\partial X_x}$ of 
the boundary $\partial X_x$ of $X_x$, and we have 
\begin{equation} \label{eq:basis1}
(\CO^x)^{\vee} = 
 \sum_{\substack{y \in W \\ y \leq x}} (-1)^{\ell(x) - \ell(y)} \CO_y, \qquad 
\CO_x = \sum_{\substack{y \in W \\ y \leq x}} (\CO^{y})^{\vee}. 
\end{equation}
Also, for $x \in W$, $(\CO_{x})^{\vee}$ is 
the class of the ideal sheaf $\CI_{\partial X^x}$ of 
the boundary $\partial X^x$ of $X^x$, and we have
\begin{equation} \label{eq:basis2}
(\CO_x)^{\vee} = \sum_{\substack{y \in W \\ y \geq x}} 
(-1)^{\ell(x) - \ell(y)} \CO^y, \qquad 
\CO^x = \sum_{\substack{y \in W \\ y \geq x}} (\CO_y)^{\vee}.
\end{equation}
Here, $\leq$ denotes the Bruhat order on $W$. 
%
%
\begin{ex}[{cf. \cite[Section~4.2]{LM}}] \label{ex:G2-B}
As in Examples~\ref{ex:G2-1} and \ref{ex:G2-2}, 
we assume that $\Fg$ is of type $G_{2}$, and $i = 2$; 
note that $\pair{\varpi_{2}}{\theta^{\vee}} \not= 1$, 
where $\theta \in \Delta^{+}$ is the highest root. 
Let $w=s_{2}s_{1}s_{2}s_{1}s_{2}$, and  
$d = d_{1}\alpha_{1}^{\vee}+2\alpha_{2}^{\vee} \in Q^{\vee,+}$ 
with $d_{1} \in \BZ_{> 0}$. Since 
\begin{equation*}
\langle \CO^{s_2},\CO^{w},(\CO^{x})^{\vee} \rangle_{d} =  
\sum_{y \le x} (-1)^{\ell(x) - \ell(y)} 
\thp{s_2}{w}{y}_{d}, 
\end{equation*}
we have
\begin{align*}
& \langle \CO^{s_2},\CO^{w},(\CO^{e})^{\vee} \rangle_{d} = 
\thp{s_2}{w}{e}_{d} = 1+\be^{-(3\alpha_{1}+2\alpha_{2})}, \\
& \langle \CO^{s_2},\CO^{w},(\CO^{s_{2}})^{\vee} \rangle_{d} = 
\thp{s_2}{w}{s_{2}}_{d} - \thp{s_2}{w}{e}_{d} = 1-(1+\be^{-\vpi_{2}}) = -\be^{-(3\alpha_{1}+2\alpha_{2})}. 
\end{align*}
Also, since 
\begin{equation*}
\thp{s_2}{w}{s_{1}}_{d} - \thp{s_2}{w}{e}_{d} = 
\twp{w}{s_{1}}_{d} - \twp{w}{e}_{d}, 
\end{equation*}
\begin{equation*}
\thp{s_2}{w}{y}_{d} = 
\twp{w}{y}_{d} \quad \text{for all $y \in W \setminus \{e,s_2\}$}, 
\end{equation*}
it follows that for $x \in W \setminus \{s_{2},e\}$ (note that $x \ge s_{1} > e$),
\begin{align*}
\langle \CO^{s_2},\CO^{w},(\CO^{x})^{\vee} \rangle_{d} & =  
\sum_{y \le x} (-1)^{\ell(x) - \ell(y)} 
\thp{s_2}{w}{y}_{d} \\
& =
\sum_{y \le x} (-1)^{\ell(x) - \ell(y)} 
\twp{w}{y}_{d} =
\sum_{y \le x} (-1)^{\ell(x) - \ell(y)} = 0.
\end{align*}

We now recall the well-known identities 
$\CO^{s_i} = 1 - \be^{- \vpi_i} \CO(- \vpi_i)$ and 
$\CO_{\lng s_i} = 1 - \be^{- \lng \vpi_i} \CO(- \vpi_i)$ in $K_{T}(X)$ for $i \in I$, 
where $\CO(- \vpi_i) \in K_{T}(X)$ is the line bundle over $X$ associated to $-\vpi_{i}$. 
Combining these, we obtain 
\begin{equation} \label{eq:line}
\CO_{\lng s_i} = (1-\be^{\vpi_i - \lng \vpi_i }) + \be^{\vpi_i - \lng \vpi_i } \CO^{s_i}. 
\end{equation}
Let $w$ and $d$ as above. From the above, 
together with equation \eqref{eq:line}, 
we compute as follows (cf. Remark~\ref{rem:positivity}): 
\begin{align*}
\langle \CO_{\lng s_2},\CO^{w},(\CO^{e})^{\vee} \rangle_{d}
 & = (1-\be^{\vpi_2 - \lng \vpi_2 }) \langle \CO^{w},(\CO^{e})^{\vee} \rangle_{d} +
    \be^{\vpi_2 - \lng \vpi_2 }\langle \CO^{s_2},\CO^{w},(\CO^{e})^{\vee} \rangle_{d} \\
 & = (1-\be^{\vpi_2 - \lng \vpi_2 }) + \be^{\vpi_2 - \lng \vpi_2 }(1+\be^{-(3\alpha_{1}+2\alpha_{2})}) \\
 & = 1 + \be^{3\alpha_{1}+2\alpha_{2}}, \\[1.5mm]
\langle \CO_{\lng s_2},\CO^{w},(\CO^{s_2})^{\vee} \rangle_{d}
 & = (1-\be^{\vpi_2 - \lng \vpi_2 }) \langle \CO^{w},(\CO^{s_2})^{\vee} \rangle_{d} +
    \be^{\vpi_2 - \lng \vpi_2 }\langle \CO^{s_2},\CO^{w},(\CO^{s_2})^{\vee} \rangle_{d} \\
 & = \be^{\vpi_2 - \lng \vpi_2 } ( -\be^{-(3\alpha_{1}+2\alpha_{2})} ) 
   = - \be^{3\alpha_{1}+2\alpha_{2}}, 
\end{align*}
and for $x \in W \setminus \{s_{2},e\}$, 
\begin{align*}
\langle \CO_{\lng s_2},\CO^{w},(\CO^{x})^{\vee} \rangle_{d}
 & = (1-\be^{\vpi_2 - \lng \vpi_2 }) \langle \CO^{w},(\CO^{x})^{\vee} \rangle_{d} +
    \be^{\vpi_2 - \lng \vpi_2 }\langle \CO^{s_2},\CO^{w},(\CO^{x})^{\vee} \rangle_{d} \\
 & = 0 + 0 = 0.
\end{align*}
\end{ex}

%

\end{document}

%% file: G2_QBG.tex
{\unitlength 0.1in%
\begin{picture}(49.0000,48.0000)(2.0000,-52.3500)%
\put(26.0000,-6.0000){\makebox(0,0){$\lng = s_1s_2s_1s_2s_1s_2 =s_2s_1s_2s_1s_2s_1$}}%
\put(16.0000,-12.0000){\makebox(0,0){$s_1s_2s_1s_2s_1$}}%
\put(36.0000,-12.0000){\makebox(0,0){$s_2s_1s_2s_1s_2$}}%
\put(36.0000,-20.0000){\makebox(0,0){$s_2s_1s_2s_1$}}%
\put(36.0000,-28.0000){\makebox(0,0){$s_2s_1s_2$}}%
\put(36.0000,-36.0000){\makebox(0,0){$s_2s_1$}}%
\put(36.0000,-44.0000){\makebox(0,0){$s_2$}}%
\put(26.0000,-50.0000){\makebox(0,0){$e$}}%
\put(16.0000,-20.0000){\makebox(0,0){$s_1s_2s_1s_2$}}%
\put(16.0000,-28.0000){\makebox(0,0){$s_1s_2s_1$}}%
\put(16.0000,-36.0000){\makebox(0,0){$s_1s_2$}}%
\put(16.0000,-44.0000){\makebox(0,0){$s_1$}}%
%
\special{pn 8}%
\special{pa 2400 4800}%
\special{pa 1800 4600}%
\special{fp}%
\special{sh 1}%
\special{pa 1800 4600}%
\special{pa 1857 4640}%
\special{pa 1851 4617}%
\special{pa 1870 4602}%
\special{pa 1800 4600}%
\special{fp}%
%
\special{pn 8}%
\special{pa 2800 4800}%
\special{pa 3400 4600}%
\special{fp}%
\special{sh 1}%
\special{pa 3400 4600}%
\special{pa 3330 4602}%
\special{pa 3349 4617}%
\special{pa 3343 4640}%
\special{pa 3400 4600}%
\special{fp}%
%
\special{pn 8}%
\special{pa 1800 4800}%
\special{pa 2400 5000}%
\special{dt 0.045}%
\special{sh 1}%
\special{pa 2400 5000}%
\special{pa 2343 4960}%
\special{pa 2349 4983}%
\special{pa 2330 4998}%
\special{pa 2400 5000}%
\special{fp}%
%
\special{pn 8}%
\special{pa 3400 4800}%
\special{pa 2800 5000}%
\special{dt 0.045}%
\special{sh 1}%
\special{pa 2800 5000}%
\special{pa 2870 4998}%
\special{pa 2851 4983}%
\special{pa 2857 4960}%
\special{pa 2800 5000}%
\special{fp}%
%
\special{pn 8}%
\special{pa 3700 4200}%
\special{pa 3700 3800}%
\special{fp}%
\special{sh 1}%
\special{pa 3700 3800}%
\special{pa 3680 3867}%
\special{pa 3700 3853}%
\special{pa 3720 3867}%
\special{pa 3700 3800}%
\special{fp}%
%
\special{pn 8}%
\special{pa 3700 3400}%
\special{pa 3700 3000}%
\special{fp}%
\special{sh 1}%
\special{pa 3700 3000}%
\special{pa 3680 3067}%
\special{pa 3700 3053}%
\special{pa 3720 3067}%
\special{pa 3700 3000}%
\special{fp}%
%
\special{pn 8}%
\special{pa 3700 2600}%
\special{pa 3700 2200}%
\special{fp}%
\special{sh 1}%
\special{pa 3700 2200}%
\special{pa 3680 2267}%
\special{pa 3700 2253}%
\special{pa 3720 2267}%
\special{pa 3700 2200}%
\special{fp}%
%
\special{pn 8}%
\special{pa 3700 1800}%
\special{pa 3700 1400}%
\special{fp}%
\special{sh 1}%
\special{pa 3700 1400}%
\special{pa 3680 1467}%
\special{pa 3700 1453}%
\special{pa 3720 1467}%
\special{pa 3700 1400}%
\special{fp}%
%
\special{pn 8}%
\special{pa 1500 1800}%
\special{pa 1500 1400}%
\special{fp}%
\special{sh 1}%
\special{pa 1500 1400}%
\special{pa 1480 1467}%
\special{pa 1500 1453}%
\special{pa 1520 1467}%
\special{pa 1500 1400}%
\special{fp}%
%
\special{pn 8}%
\special{pa 1500 2600}%
\special{pa 1500 2200}%
\special{fp}%
\special{sh 1}%
\special{pa 1500 2200}%
\special{pa 1480 2267}%
\special{pa 1500 2253}%
\special{pa 1520 2267}%
\special{pa 1500 2200}%
\special{fp}%
%
\special{pn 8}%
\special{pa 1500 3400}%
\special{pa 1500 3000}%
\special{fp}%
\special{sh 1}%
\special{pa 1500 3000}%
\special{pa 1480 3067}%
\special{pa 1500 3053}%
\special{pa 1520 3067}%
\special{pa 1500 3000}%
\special{fp}%
%
\special{pn 8}%
\special{pa 1500 4200}%
\special{pa 1500 3800}%
\special{fp}%
\special{sh 1}%
\special{pa 1500 3800}%
\special{pa 1480 3867}%
\special{pa 1500 3853}%
\special{pa 1520 3867}%
\special{pa 1500 3800}%
\special{fp}%
%
\special{pn 8}%
\special{pa 3500 3800}%
\special{pa 3500 4200}%
\special{dt 0.045}%
\special{sh 1}%
\special{pa 3500 4200}%
\special{pa 3520 4133}%
\special{pa 3500 4147}%
\special{pa 3480 4133}%
\special{pa 3500 4200}%
\special{fp}%
%
\special{pn 8}%
\special{pa 3500 3000}%
\special{pa 3500 3400}%
\special{dt 0.045}%
\special{sh 1}%
\special{pa 3500 3400}%
\special{pa 3520 3333}%
\special{pa 3500 3347}%
\special{pa 3480 3333}%
\special{pa 3500 3400}%
\special{fp}%
%
\special{pn 8}%
\special{pa 3500 2200}%
\special{pa 3500 2600}%
\special{dt 0.045}%
\special{sh 1}%
\special{pa 3500 2600}%
\special{pa 3520 2533}%
\special{pa 3500 2547}%
\special{pa 3480 2533}%
\special{pa 3500 2600}%
\special{fp}%
%
\special{pn 8}%
\special{pa 3500 1400}%
\special{pa 3500 1800}%
\special{dt 0.045}%
\special{sh 1}%
\special{pa 3500 1800}%
\special{pa 3520 1733}%
\special{pa 3500 1747}%
\special{pa 3480 1733}%
\special{pa 3500 1800}%
\special{fp}%
%
\special{pn 8}%
\special{pa 1700 1400}%
\special{pa 1700 1800}%
\special{dt 0.045}%
\special{sh 1}%
\special{pa 1700 1800}%
\special{pa 1720 1733}%
\special{pa 1700 1747}%
\special{pa 1680 1733}%
\special{pa 1700 1800}%
\special{fp}%
%
\special{pn 8}%
\special{pa 1700 2200}%
\special{pa 1700 2600}%
\special{dt 0.045}%
\special{sh 1}%
\special{pa 1700 2600}%
\special{pa 1720 2533}%
\special{pa 1700 2547}%
\special{pa 1680 2533}%
\special{pa 1700 2600}%
\special{fp}%
%
\special{pn 8}%
\special{pa 1700 3000}%
\special{pa 1700 3400}%
\special{dt 0.045}%
\special{sh 1}%
\special{pa 1700 3400}%
\special{pa 1720 3333}%
\special{pa 1700 3347}%
\special{pa 1680 3333}%
\special{pa 1700 3400}%
\special{fp}%
%
\special{pn 8}%
\special{pa 1700 3800}%
\special{pa 1700 4200}%
\special{dt 0.045}%
\special{sh 1}%
\special{pa 1700 4200}%
\special{pa 1720 4133}%
\special{pa 1700 4147}%
\special{pa 1680 4133}%
\special{pa 1700 4200}%
\special{fp}%
%
\special{pn 8}%
\special{pa 3700 1000}%
\special{pa 3100 800}%
\special{fp}%
\special{sh 1}%
\special{pa 3100 800}%
\special{pa 3157 840}%
\special{pa 3151 817}%
\special{pa 3170 802}%
\special{pa 3100 800}%
\special{fp}%
%
\special{pn 8}%
\special{pa 2800 800}%
\special{pa 3400 1000}%
\special{dt 0.045}%
\special{sh 1}%
\special{pa 3400 1000}%
\special{pa 3343 960}%
\special{pa 3349 983}%
\special{pa 3330 998}%
\special{pa 3400 1000}%
\special{fp}%
%
\special{pn 8}%
\special{pa 1500 1000}%
\special{pa 2100 800}%
\special{fp}%
\special{sh 1}%
\special{pa 2100 800}%
\special{pa 2030 802}%
\special{pa 2049 817}%
\special{pa 2043 840}%
\special{pa 2100 800}%
\special{fp}%
%
\special{pn 8}%
\special{pa 2500 800}%
\special{pa 1900 1000}%
\special{dt 0.045}%
\special{sh 1}%
\special{pa 1900 1000}%
\special{pa 1970 998}%
\special{pa 1951 983}%
\special{pa 1957 960}%
\special{pa 1900 1000}%
\special{fp}%
\put(28.0000,-43.5000){\makebox(0,0)[lb]{$\clubsuit$}}%
\put(28.0000,-35.5000){\makebox(0,0)[lb]{$\diamondsuit$}}%
%
\special{pn 8}%
\special{pa 1920 4400}%
\special{pa 3320 3700}%
\special{fp}%
\special{sh 1}%
\special{pa 3320 3700}%
\special{pa 3251 3712}%
\special{pa 3272 3724}%
\special{pa 3269 3748}%
\special{pa 3320 3700}%
\special{fp}%
%
\special{pn 8}%
\special{pa 3300 4400}%
\special{pa 1900 3700}%
\special{fp}%
\special{sh 1}%
\special{pa 1900 3700}%
\special{pa 1951 3748}%
\special{pa 1948 3724}%
\special{pa 1969 3712}%
\special{pa 1900 3700}%
\special{fp}%
%
\special{pn 8}%
\special{pa 3300 3600}%
\special{pa 1900 2900}%
\special{fp}%
\special{sh 1}%
\special{pa 1900 2900}%
\special{pa 1951 2948}%
\special{pa 1948 2924}%
\special{pa 1969 2912}%
\special{pa 1900 2900}%
\special{fp}%
%
\special{pn 8}%
\special{pa 3300 2800}%
\special{pa 1900 2100}%
\special{fp}%
\special{sh 1}%
\special{pa 1900 2100}%
\special{pa 1951 2148}%
\special{pa 1948 2124}%
\special{pa 1969 2112}%
\special{pa 1900 2100}%
\special{fp}%
%
\special{pn 8}%
\special{pa 3300 2000}%
\special{pa 1900 1300}%
\special{fp}%
\special{sh 1}%
\special{pa 1900 1300}%
\special{pa 1951 1348}%
\special{pa 1948 1324}%
\special{pa 1969 1312}%
\special{pa 1900 1300}%
\special{fp}%
%
\special{pn 8}%
\special{pa 1920 3600}%
\special{pa 3320 2900}%
\special{fp}%
\special{sh 1}%
\special{pa 3320 2900}%
\special{pa 3251 2912}%
\special{pa 3272 2924}%
\special{pa 3269 2948}%
\special{pa 3320 2900}%
\special{fp}%
%
\special{pn 8}%
\special{pa 1920 2800}%
\special{pa 3320 2100}%
\special{fp}%
\special{sh 1}%
\special{pa 3320 2100}%
\special{pa 3251 2112}%
\special{pa 3272 2124}%
\special{pa 3269 2148}%
\special{pa 3320 2100}%
\special{fp}%
%
\special{pn 8}%
\special{pa 1920 2000}%
\special{pa 3320 1300}%
\special{fp}%
\special{sh 1}%
\special{pa 3320 1300}%
\special{pa 3251 1312}%
\special{pa 3272 1324}%
\special{pa 3269 1348}%
\special{pa 3320 1300}%
\special{fp}%
\put(20.0000,-33.5000){\makebox(0,0)[lt]{$\spadesuit$}}%
\put(20.0000,-25.5000){\makebox(0,0)[lt]{$\spadesuit$}}%
\put(28.0000,-27.5000){\makebox(0,0)[lb]{$\diamondsuit$}}%
\put(28.0000,-19.5000){\makebox(0,0)[lb]{$\clubsuit$}}%
%
\special{pn 8}%
\special{pa 1200 2800}%
\special{pa 1000 2800}%
\special{dt 0.045}%
%
\special{pn 8}%
\special{pa 1000 2800}%
\special{pa 1000 5200}%
\special{dt 0.045}%
%
\special{pn 8}%
\special{pa 1000 1200}%
\special{pa 600 1200}%
\special{dt 0.045}%
%
\special{pn 8}%
\special{pa 600 1200}%
\special{pa 600 3600}%
\special{dt 0.045}%
%
\special{pn 8}%
\special{pa 600 3600}%
\special{pa 1200 3600}%
\special{dt 0.045}%
\special{sh 1}%
\special{pa 1200 3600}%
\special{pa 1133 3580}%
\special{pa 1147 3600}%
\special{pa 1133 3620}%
\special{pa 1200 3600}%
\special{fp}%
%
\special{pn 8}%
\special{pa 4100 2000}%
\special{pa 4300 2000}%
\special{dt 0.045}%
%
\special{pn 8}%
\special{pa 4300 2000}%
\special{pa 4300 4400}%
\special{dt 0.045}%
%
\special{pn 8}%
\special{pa 4300 4400}%
\special{pa 3700 4400}%
\special{dt 0.045}%
\special{sh 1}%
\special{pa 3700 4400}%
\special{pa 3767 4420}%
\special{pa 3753 4400}%
\special{pa 3767 4380}%
\special{pa 3700 4400}%
\special{fp}%
%
\special{pn 8}%
\special{pa 4000 600}%
\special{pa 4800 600}%
\special{dt 0.045}%
%
\special{pn 8}%
\special{pa 4800 600}%
\special{pa 4800 2800}%
\special{dt 0.045}%
%
\special{pn 8}%
\special{pa 4800 2800}%
\special{pa 4000 2800}%
\special{dt 0.045}%
\special{sh 1}%
\special{pa 4000 2800}%
\special{pa 4067 2820}%
\special{pa 4053 2800}%
\special{pa 4067 2780}%
\special{pa 4000 2800}%
\special{fp}%
%
\special{pn 8}%
\special{pa 1000 600}%
\special{pa 200 600}%
\special{dt 0.045}%
%
\special{pn 8}%
\special{pa 200 600}%
\special{pa 200 4400}%
\special{dt 0.045}%
%
\special{pn 8}%
\special{pa 200 4400}%
\special{pa 1400 4400}%
\special{dt 0.045}%
\special{sh 1}%
\special{pa 1400 4400}%
\special{pa 1333 4380}%
\special{pa 1347 4400}%
\special{pa 1333 4420}%
\special{pa 1400 4400}%
\special{fp}%
\put(6.0000,-5.0000){\makebox(0,0){$\spadesuit$}}%
%
\special{pn 8}%
\special{pa 1000 5200}%
\special{pa 2500 5200}%
\special{dt 0.045}%
%
\special{pn 8}%
\special{pa 2500 5200}%
\special{pa 2500 5100}%
\special{dt 0.045}%
\special{sh 1}%
\special{pa 2500 5100}%
\special{pa 2480 5167}%
\special{pa 2500 5153}%
\special{pa 2520 5167}%
\special{pa 2500 5100}%
\special{fp}%
%
\special{pn 8}%
\special{pa 2700 5200}%
\special{pa 2700 5100}%
\special{dt 0.045}%
\special{sh 1}%
\special{pa 2700 5100}%
\special{pa 2680 5167}%
\special{pa 2700 5153}%
\special{pa 2720 5167}%
\special{pa 2700 5100}%
\special{fp}%
%
\special{pn 8}%
\special{pa 4100 1200}%
\special{pa 5100 1200}%
\special{dt 0.045}%
%
\special{pn 8}%
\special{pa 5100 1200}%
\special{pa 5100 5200}%
\special{dt 0.045}%
%
\special{pn 8}%
\special{pa 5100 5200}%
\special{pa 2700 5200}%
\special{dt 0.045}%
\put(44.0000,-53.0000){\makebox(0,0){$\spadesuit$}}%
\end{picture}}%